\documentclass[a4paper,12pt]{article} 

\usepackage{color}

\usepackage{amsfonts, amsmath, amsthm, amssymb}
\usepackage[T1]{fontenc}
\usepackage[cp1250]{inputenc}
\usepackage{xcolor}
\usepackage{graphicx}
\usepackage{amssymb}
\usepackage{amsmath}
\usepackage{mathptmx}
\usepackage{helvet}
\usepackage{courier}
\usepackage{txfonts}
\usepackage{tikz} 
\usetikzlibrary{arrows}
\usepackage{type1cm}

\usepackage{verbatim}

\usepackage{graphicx}
\usepackage{epsfig,amscd,amssymb,amsxtra,amsmath,amsthm}
\usepackage{type1cm}
\usepackage[T1]{fontenc}
\usepackage{graphics}
\usepackage[mathscr]{eucal}
\usepackage[all]{xy}
\usepackage{amsmath,amscd}

\newcommand{\orbit}{\mathcal{O}_f^{\oplus}}
\newcommand{\orbitminus}{\mathcal{O}_f^{\ominus}}
\newcommand{\Orbitminus}{\mathcal{O}_G^{\ominus}}
\newcommand{\Orbit}{\mathcal{O}_G^{\oplus}}
\newcommand{\Orbitt}{\mathcal{U}_G^{\oplus}}

\newtheorem{theorem}{Theorem}[section]
\newtheorem{proposition}[theorem]{Proposition}
\newtheorem{definition}[theorem]{Definition}
\newtheorem{lemma}[theorem]{Lemma}

\newtheorem{example}[theorem]{Example}

\newtheorem{problem}[theorem]{Problem}
\newtheorem{observation}[theorem]{Observation}

\newcommand{\Cl}  {\mathop{\rm Cl}\nolimits}

\begin{document}

\def\joinrel{\mkern-3mu}
\newcommand{\varproj}{\displaystyle \lim_{\multimapinv\joinrel-\joinrel-}}

\title{Minimal dynamical systems with closed relations}
\author{Iztok Bani\v c,  Goran Erceg, Rene Gril Rogina and Judy Kennedy}
\date{}

\maketitle

\begin{abstract}
\noindent We introduce dynamical systems $(X,G)$ with closed relations $G$ on compact metric spaces $X$ and discuss different types of  minimality of such dynamical systems, all of them generalizing minimal dynamical systems $(X,f)$ with continuous function $f$ on a compact metric space $X$.
\end{abstract}
\-


\section{Introduction}\label{s000}

Our work is motivated by E.~Akin's book ``General Topology of Dynamical Systems'' \cite{A}, where dynamical systems using closed relations are presented, and by S.~Kolyada's and L.~Snoha's paper  ``Minimal dynamical systems'' \cite{KS}, where a wonderful overview of minimal dynamical systems is given.  Our work is also motivated by  J.~Kennedy's and G.~Erceg's paper ``Topological entropy on closed sets in $[0,1]^2$'' \cite{EK}, and I.~Bani\v c's,  J.~Kennedy's and G.~Erceg's paper ``Closed relations with non-zero entropy that generate no periodic points'' \cite{BEK}, where the idea of topological entropy is generalized from standard dynamical systems $(X,f)$ to dynamical systems $(X,G)$ with closed relations $G$ on compact metric spaces $X$.

 In dynamical systems theory,  the study of chaotic behaviour of a dynamical system is  often based on some topological properties or properties of continuous functions. One of the commonly used properties is the minimality of a dynamical system $(X,f)$ or the minimality of the function $f$. According to \cite{KS}, minimal dynamical systems were defined by Birkhoff in 1912 \cite{B} as the systems which have no nontrivial closed subsystems: they are considered to be the most fundamental dynamical systems;   see \cite{KS} where more references can be found.  Minimal dynamical systems $(X,f)$ (i.e., with a minimal map $f$) have the property that each point moves under iteration of $f$ from one non-empty open set to another. This property has been studied intensively by mathematicians since it is an important property in dynamical system theory.

In this paper, we generalize the notion of topological dynamical systems to topological dynamical systems with closed relations and introduce the notion of minimality of such dynamical systems.   A similar generalization of a dynamical object was presented in 2004 by Ingram and Mahavier  \cite{ingram,mah} introducing inverse limits of inverse sequences of compact metric spaces $X$ with upper semi-continuous set-valued bonding functions $f$ (their graphs $\Gamma(f)$ are examples of closed relations on $X$ with certain additional properties).  These inverse limits provide a valuable extension to the role of inverse limits in the study of dynamical systems and continuum theory.  For example,   Kennedy and Nall  have developed a simple method for constructing families of $\lambda$-dendroids \cite{KN}. Their method involves inverse limits of inverse sequences with upper semi-continuous set-valued functions on closed intervals with simple bonding functions. 
Such generalizations have proven to be  useful (also in applied areas); frequently,  when  constructing a model for  empirical data, continuous (single-valued) functions fall short, and the data are better modelled by upper semi-continuous set-valued functions, or sometimes, even closed relations that are not set-valued functions are required. The Christiano-Harrison model   from macroeconomics is one such example \cite{christiano}.  The study of inverse limits of inverse sequences with upper semi-continuous set-valued  functions is rapidly gaining momentum - the recent books by Ingram \cite{ingram-knjiga}, and by Ingram and Mahavier \cite{ingram}, give a comprehensive exposition of this research prior to 2012.   

Also,  several papers on the topic of dynamical systems with (upper semi-continuous) set-valued functions have appeared recently, see \cite{CP,LP,LYY,LWZ,KN,KW,MRT,R,SWS}, where more references may be found.  However,  there is not much known of such dynamical systems and therefore,  there are many properties of such set-valued dynamical systems that are yet to be studied.  In this paper, we study the minimality of such dynamical systems.  We also extend the notion of dynamical systems with (upper semi-continuous) set-valued functions to dynamical systems with closed relations.

We proceed as follows. In the sections that follow Section \ref{s00}, where basic definitions are given, we discuss the following topics:
\begin{enumerate}
\item Minimal dynamical systems with closed relations and invariant sets (Section \ref{s1}).
\item Minimal dynamical systems with closed relations and forward orbits (Section \ref{s2}).
\item Minimal dynamical systems with closed relations and omega limit sets (Section \ref{s3}).
\item Backward minimal dynamical systems with closed relations (Section \ref{s4}).
\item Minimal dynamical systems with closed relations and backward orbits (Section \ref{s5}).
\item Minimal dynamical systems with closed relations and alpha limit sets (Section \ref{s6}).
\item Preserving minimality by topological conjugation (Section \ref{s8}).
\end{enumerate}
In  Sections \ref{s1},  \ref{s2},  \ref{s3},  \ref{s5},  \ref{s6},   and \ref{s8},   we first revisit minimal dynamical systems $(X,f)$ and then, we  generalize the asserted property from dynamical systems $(X,f)$ to dynamical systems with closed relations $(X,G)$ by making the identification $(X,f)=(X,\Gamma(f))$.   Results about dynamical systems $(X,f)$,  presented in first part of each of the above mentioned sections, are well-known. Their proofs are short, rather straight forward and elementary.  Since they are important for the purpose of this paper,  we state and prove each of the presented results.

\section{Definitions and notation} \label{s00}

In this section, basic definitions and well-known results that are needed later in the paper are presented.
\begin{definition}
Let $X$ and $Y$ be metric spaces, and let $f:X\rightarrow Y$ be a function.  We use  
$$
\Gamma(f)=\{(x,y)\in X\times Y \ | \ y=f(x)\}
$$
to denote \emph{ \color{blue} the graph of the function $f$}.
\end{definition}
\begin{definition}
If $X$ is a compact metric space, then \emph{ \color{blue}  $2^X$ }denotes the set of all  non-empty closed subsets of $X$.  
\end{definition}
%

\begin{definition}
Let $X$ be a compact metric space and let $G\subseteq X\times X$ be a relation on $X$. If $G\in 2^{X\times X}$, then we say that $G$ is \emph{ \color{blue} a closed relation on $X$}.  
\end{definition}

\begin{definition}
Let $X$  be a set and let $G$ be a relation on $X$.  Then we define  
$$
G^{-1}=\{(y,x)\in X\times X \ | \ (x,y)\in G\}
$$
to be \emph{ \color{blue} the inverse relation of the relation $G$} on $X$.
\end{definition}
\begin{definition}
Let $X$ be a compact metric space and let $G$ be a closed relation on $X$. Then we call
$$
\star_{i=1}^{m}G=\Big\{(x_1,x_2,x_3,\ldots ,x_{m+1})\in \prod_{i=1}^{m+1}X \ | \ \textup{ for each } i\in \{1,2,3,\ldots ,m\}, (x_{i},x_{i+1})\in G\Big\}
$$
for each positive integer $m$, \emph{ \color{blue} the $m$-th Mahavier product of $G$}, and
$$
\star_{i=1}^{\infty}G=\Big\{(x_1,x_2,x_3,\ldots )\in \prod_{i=1}^{\infty}X \ | \ \textup{ for each positive integer } i, (x_{i},x_{i+1})\in G\Big\}
$$
\emph{ \color{blue} the infinite  Mahavier product of $G$}.
\end{definition}
\begin{observation}
Let $X$ be a compact metric space, let $f:X\rightarrow X$ be a continuous function. 
Then 
$$
\star_{n=1}^{\infty}\Gamma(f)^{-1}=\varprojlim(X,f).
$$
\end{observation}

\section{Minimal dynamical systems with closed relations}\label{s1}
First, we revisit minimal dynamical systems and then, we introduce dynamical systems with closed relations and generalize the notion of minimality of a dynamical system to minimality of dynamical systems with closed relations.  
\begin{definition}
Let $X$ be a compact metric space, let $f:X\rightarrow X$ be a continuous function. We say that $(X,f)$ is \emph{ \color{blue} a dynamical system}.
\end{definition}
\begin{definition}
Let $(X,f)$ be a dynamical system and let $A\subseteq X$. We say that 
\begin{enumerate}
\item $A$ is \emph{ \color{blue} $f$-invariant}, if $f(A)\subseteq A$.
\item $A$ is \emph{ \color{blue} strongly $f$-invariant}, if $f(A)=A$.
\end{enumerate}
\end{definition}

\begin{definition}
Let $(X,f)$ be a dynamical system.   We  say that $(X,f)$ is \emph{ \color{blue} a minimal dynamical system}, if for each closed subset $A$ of  $X$, 
$$
A \textup{ is } f\textup{-invariant}  \Longrightarrow A\in \{\emptyset ,X\}.
$$
\end{definition}

The following are well-known results. Since the proofs are short and elementary, we give them here.
\begin{theorem}\label{enolicni}
Let $(X,f)$ be a dynamical system.   The following statements are equivalent.
\begin{enumerate}
\item \label{ena} $(X,f)$ is a minimal dynamical system.
\item\label{dva} For each closed subset $A$ of  $X$, 
$$
f(A)= A  \Longrightarrow A\in \{\emptyset ,X\}.
$$
\end{enumerate} 
\end{theorem}
\begin{proof}
Let $(X,f)$ be a minimal dynamical system and let $A$ be a closed subset of $X$ such that $f(A)=A$. Then $f(A)\subseteq A$. Since $(X,f)$ is minimal, it follows that $A\in \{\emptyset ,X\}$. Next, suppose that for each closed subset $A$ of  $X$, 
$$
f(A)= A  \Longrightarrow A\in \{\emptyset ,X\}.
$$
To prove that $(X,f)$ is minimal, let $A$ be a closed subset of $X$ such that $f(A)\subseteq A$.  Set
$$
\mathcal U=\{B\in 2^{A} \ | \ f(B)\subseteq B\}.
$$
We use Zorn's lemma to show that there is a minimal element of $(\mathcal U,\subseteq)$. Let $\mathcal V$ be a chain (a totally ordered set) in $\mathcal U$ and set $B_0=\bigcap_{B\in \mathcal V}B$. Then $B_0\in 2^A$ and
$$
f(B_0)=f\Big(\bigcap_{B\in \mathcal V}B\Big)\subseteq \bigcap_{B\in \mathcal V}f(B)\subseteq \bigcap_{B\in \mathcal V}B=B_0.
$$
Therefore, $B_0\in \mathcal U$  and $B_0\subseteq B$ for any $B\in \mathcal V$. By Zorn's lemma, there is a minimal element of $(\mathcal U,\subseteq)$. Let $A_0$ be a minimal element of $(\mathcal U,\subseteq)$. Then $A_0$ is a non-empty closed subset of $X$. Note that $f(A_0)\subseteq A_0$ since $A_0$ is an element of $(\mathcal U,\subseteq)$. Suppose that $f(A_0)\neq A_0$ and let $A_1=f(A_0)$. Then $A_1$ is a proper subset of $A_0$ such that $f(A_1)\subseteq A_1$ (note that $f(A_1)=f(f(A_0))\subseteq f(A_0)= A_1$ since $f(A_0)\subseteq A_0$). This contradicts the minimality of $A_0$ in  $(\mathcal U,\subseteq)$. Therefore, $A_0$ is a non-empty closed subset of $X$ such that $f(A_0)=A_0$.  It follows that $A_0=X$ since $A_0\neq \emptyset$.  Therefore,  $A=X$ since $A_0\subseteq A$. We have just proved that  \ref{ena} is equivalent to \ref{dva}. 
\end{proof}

Next, we introduce dynamical systems with closed relations.  Before we do that, we give an obvious Proposition \ref{prop1cc}, which will serve as a motivation for the rest of this section.

\begin{definition}
Let $X$ be a metric space.  We  use 
$$
p_1,p_2:X\times X\rightarrow X
$$
to denote\emph{ \color{blue}  the standard projections} defined by
$$
p_1(x,y)=x  \textup{ and } p_2(x,y)=y
$$
for all $(x,y)\in X\times X$.
\end{definition}
\begin{proposition}\label{prop1cc}
Let $(X,f)$ be a dynamical system and let $A\subseteq X$.
The following statements are equivalent:
\begin{enumerate}
\item $A$ is $f$-invariant.
\item For each $(x,y)\in \Gamma(f)$,
$$
x\in A \Longrightarrow y\in A.
$$
\item For each $x\in A$,
$$
x\in p_1(\Gamma(f)) \Longrightarrow \textup{ there is  } y\in A \textup{ such that } (x,y)\in \Gamma(f).
$$
\end{enumerate}
\end{proposition}
\begin{proof}
Suppose that $A$ is $f$ invariant and let $(x,y)\in \Gamma(f)$ such that $x\in A$. Then $y=f(x)$ and since $A$ is $f$ invariant, it follows that $y\in A$. 

Suppose that for each $(x,y)\in \Gamma(f)$,
$$
x\in A \Longrightarrow y\in A
$$
and let $x\in A$ such that $x\in p_1(\Gamma(f))$. Set $y=f(x)$. Then $(x,y)\in \Gamma(f)$ and $y\in A$ follows.

Suppose that for each $x\in A$,
$$
x\in p_1(\Gamma(f)) \Longrightarrow \textup{ there is  } y\in A \textup{ such that } (x,y)\in \Gamma(f)
$$
and let $x\in A$. Let $y\in A$ such that $(x,y)\in \Gamma(f)$. Then $y=f(x)$ and $f(x)\in A$ follows. Therefore, $A$ is $f$-invariant.
\end{proof}

Motivated by Proposition \ref{prop1cc},  we  introduce two different types of invariant sets with respect to a closed relation on a compact metric space. 
\begin{definition}
Let $X$ be a compact metric space and let $G$ be a non-empty closed relation on $X$. We say that $(X,G)$ is \emph{ \color{blue}  a dynamical system with a closed relation} or, briefly,  \emph{ \color{blue}  a CR-dynamical system}. 
\end{definition}

\begin{definition}
Let $(X,G)$ be a CR-dynamical system and let $A\subseteq X$. We say that the set $A$ is 
\begin{enumerate}
\item \emph{ \color{blue}  $1$-invariant in $(X,G)$},  if for each $x\in A$,
$$
x\in p_1(G) \Longrightarrow \textup{ there is  } y\in A \textup{ such that } (x,y)\in G.
$$
\item \emph{ \color{blue}  $\infty$-invariant in $(X,G)$}, if for each $(x,y)\in G$,
$$
x\in A \Longrightarrow y\in A.
$$
\end{enumerate} 
\end{definition}

\begin{observation}\label{gorane1}
Let $(X,G)$ be a CR-dynamical system,  let $\mathbf x=(x_1,x_2,x_3,\ldots ) \in \star_{i=1}^{\infty}G$, and let  let $A$ be an $\infty$-invariant set in $(X,G)$.  If $x_1\in A$, then $x_k\in A$ for any positive integer $k$.
\end{observation}

\begin{observation}\label{obs1}
Let $(X,f)$ be a dynamical system and let $A\subseteq X$. Then $(X,\Gamma(f))$ is a CR-dynamical system and by Proposition \ref{prop1cc}, the following statements are equivalent.
\begin{enumerate}
\item The set $A$ is $f$-invariant.
\item The set $A$ is $1$-invariant in $(X,\Gamma(f))$.
\item The set $A$ is $\infty$-invariant in $(X,\Gamma(f))$.
\end{enumerate}
\end{observation}
Next, we show that every $\infty$-invariant in $(X,G)$ set is also $1$-invariant in $(X,G)$.
\begin{proposition}\label{mimika}
Let $(X,G)$ be a CR-dynamical system and let $A\subseteq X$. If $A$ is $\infty$-invariant in $(X,G)$, then $A$ is $1$-invariant in $(X,G)$.
\end{proposition}
\begin{proof}
Suppose that $A$ is $\infty$-invariant in $(X,G)$. If $A\cap p_1(G)=\emptyset$, then there is nothing to show, so, let $x\in A\cap p_1(G)$ and let $y\in X$ be any point such that $(x,y)\in G$. Such a point exists since $x\in p_1(G)$. Since $x\in A$ and since $A$ is $\infty$-invariant in $(X,G)$, $y\in A$. So, there is a point $y\in A$ such that $(x,y)\in G$. Therefore, $A$ is $1$-invariant in $(X,G)$.
\end{proof}
The following example shows that there are CR-dynamical systems $(X,G)$ and subsets $A$ of $X$ such that $A$ is $1$-invariant in $(X,G)$ but it is not $\infty$-invariant in $(X,G)$.
\begin{figure}[h!]
	\centering
		\includegraphics[width=13em]{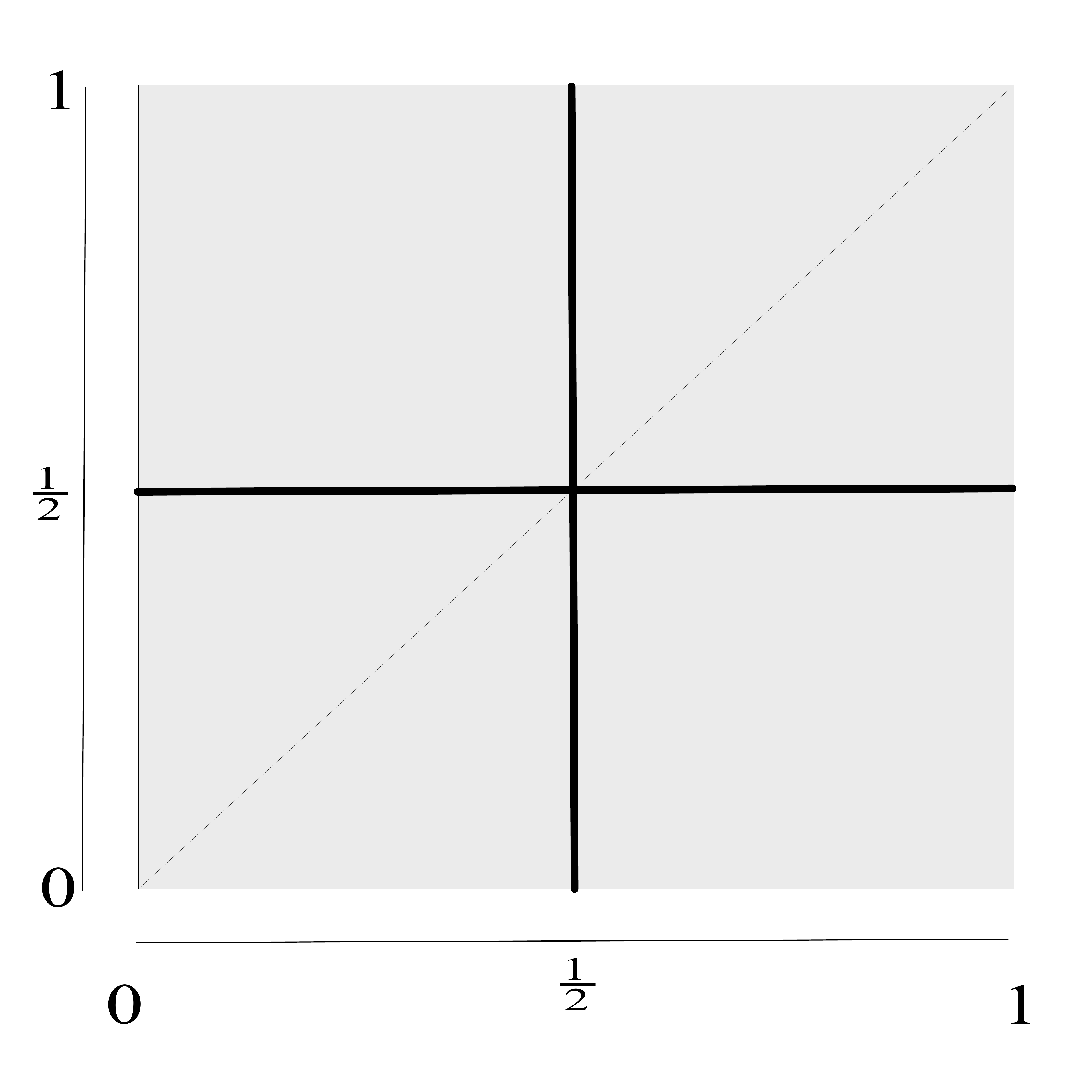}
	\caption{The relation  $G$ from Example \ref{ex1}}
	\label{figure1}
\end{figure} 
\begin{example} \label{ex1}
Let $X=[0,1]$ and let $G=([0,1]\times \{\frac{1}{2}\})\cup(\{\frac{1}{2}\}\times [0,1])$, see Figure \ref{figure1}.
Then $(X,G)$ is a CR-dynamical system. Let $A=\{\frac{1}{2}\}$. Then $A$ is $1$-invariant in $(X,G)$ but it is not $\infty$-invariant in $(X,G)$.
\end{example}
 
\begin{definition}
Let $(X,G)$ be a CR-dynamical system. We say that  
\begin{enumerate}
\item $(X,G)$ is  \emph{ \color{blue}  $1$-minimal} if for each closed subset $A$ of $X$,
$$
A \textup{ is } 1\textup{-invariant in } (X,G) \Longrightarrow A\in \{\emptyset, X\}.
$$
\item $(X,G)$ is \emph{ \color{blue}  $\infty$-minimal} if for each closed subset $A$ of $X$,
$$
A \textup{ is } \infty\textup{-invariant in } (X,G) \Longrightarrow A\in \{\emptyset, X\}.
$$
\end{enumerate}
\end{definition}
\begin{theorem}\label{judyK}
Let $(X,G)$ be a CR-dynamical system.  If $(X,G)$ is $1$-minimal,  then $(X,G)$ is $\infty$-minimal. 
\end{theorem}
\begin{proof}
Let $(X,G)$ be $1$-minimal and let $A$ be a non-empty closed subset of $X$ such that $A$ is $\infty$-invariant in $(X,G)$.  Then $A$ is $1$-invariant in $(X,G)$ by Proposition \ref{mimika}. Therefore, $A=X$ and it follows that $(X,G)$ is $\infty$-minimal.
\end{proof}
In the following example,  we show that there is a $\infty$-minimal CR-dynamical system which is not $1$-minimal.
\begin{example}\label{goranH}
Let $X=[0,1]$ and let $G=([0,1]\times \{\frac{1}{2}\})\cup(\{\frac{1}{2}\}\times [0,1])$, see Figure \ref{figure1}. Then $(X,G)$ is $\infty$-minimal but it is not $1$-minimal. Let $A=\{\frac{1}{2}\}$.  Then $A$ is $1$-invariant in $(X,G)$ but $A\not \in \{\emptyset,X\}$. Therefore, $(X,G)$ is not $1$-minimal.  To see that $(X,G)$ is $\infty$-minimal, let $A$ be a non-empty closed subset of $X$ such that $A$ is $\infty$-invariant in $(X,G)$.  Let $x\in A$. Then $(x,\frac{1}{2})\in G$ and $\frac{1}{2}\in A$ follows.  Since $(\frac{1}{2},t)\in G$, it follows that $t\in A$ for any $t\in X$. Therefore, $A=X$ and it follows that $(X,G)$ is $\infty$-minimal.
\end{example}

\section{Minimality and forward orbits}\label{s2}
First, we revisit orbits of dynamical systems $(X,f)$ and then we generalize these to orbits of CR-dynamical systems $(X,G)$.
\begin{definition}
Let $(X,f)$ be a dynamical system and let $x_0\in X$. The sequence 
$$
\mathbf x=(x_0,f(x_0),f^2(x_0),f^3(x_0),\ldots)\in \star_{i=1}^{\infty}\Gamma(f)
$$   
is called \emph{ \color{blue} the forward orbit of $x_0$.} The set 
$$
\mathcal O_f^{\oplus}(\mathbf x)=\{x_0,f(x_0),f^2(x_0),f^3(x_0),\ldots\}
$$   
is called \emph{ \color{blue} the forward orbit set of $x_0$}.
\end{definition}

The following are well-known results. Since the proof is short and elementary, we give it here.
\begin{theorem}\label{enolicnicc}
Let $(X,f)$ be a dynamical system.   The following statements are equivalent.
\begin{enumerate}
\item \label{enacc} $(X,f)$ is a minimal dynamical system.
\item \label{tricc} For each $x\in X$, 
$$
\Cl (\orbit(x))=X.
$$
\end{enumerate} 
\end{theorem}
\begin{proof}
First, we prove that  \ref{enacc} implies \ref{tricc}. Let $(X,f)$ be a minimal dynamical system and let $x\in X$ be any point.  Let $A=\Cl(\orbit(x))$. Then $A$ is a non-empty closed subset of $X$ such that $f(A)\subseteq A$. Since $(X,f)$ is minimal, it follows that $A=X$. Therefore,  $\Cl(\orbit(x))=X$.  Finally,  suppose that for each $x\in X$,  $\Cl (\orbit(x))=X$. To show that $(X,f)$ is minimal, let $A$ be a non-empty closed subset of $X$ such that $f(A)\subseteq A$.  Let $x\in A$.  Then $\Cl(\orbit(x))\subseteq A$, and $A=X$ follows from $\Cl(\orbit(x))=X$. Therefore,  $(X,f)$ is minimal. We have proved that \ref{enacc} is equivalent to \ref{tricc}. 
\end{proof}

\begin{definition}
Let $X$ be a compact metric space.  For each positive integer $k$,  we use $\pi_k:\prod_{i=1}^{\infty}X\rightarrow X$ to denote \emph{ \color{blue}  the $k$-th standard projection} from $\prod_{i=1}^{\infty}X$ to $X$.
\end{definition}
\begin{definition}
Let $(X,G)$ be a CR-dynamical system and let $x_0\in X$. We use \emph{ \color{blue}  $T_G^{+}(x_0)$} to denote the set 
$$
T_G^{+}(x_0)=\{\mathbf x\in \star_{i=1}^{\infty}G \ | \  \pi_1(\mathbf x)=x_0\}\subseteq \star_{i=1}^{\infty}G.
$$
\end{definition}
\begin{definition}
Let $(X,G)$ be a CR-dynamical system, let $\mathbf x\in \star_{i=1}^{\infty}G$, and let $x_0\in X$. We 
 \begin{enumerate}
 \item say that $\mathbf x$ is \emph{\color{blue}a forward  orbit of $x_0$ in $(X,G)$}, if $\pi_1(\mathbf x)=x_0$. 
 \item use \emph{ \color{blue}  $\mathcal O_G^{\oplus}(\mathbf x)$} to denote the set
 $$
 \mathcal O_G^{\oplus}(\mathbf x)=\{\pi_k(\mathbf x) \ | \  k \textup{ is a positive integer}\}\subseteq X.
 $$
 \item use \emph{ \color{blue}  $\mathcal U_G^{\oplus}(x_0)$} to denote the set
 $$
 \mathcal U_G^{\oplus}(x_0)=\bigcup_{\mathbf x\in T_G^{+}(x_0)}\mathcal O_G^{\oplus}(\mathbf x)\subseteq X.
 $$
 \end{enumerate}
\end{definition}
\begin{example}
Let $X=[0,1]$ and let $G=\{(1,0)\}$. Then $\star_{i=1}^{1}G\neq \emptyset $ and for each $m\neq 1$,  $\star_{i=1}^{m}G=\emptyset$.  Therefore, in this CR-dynamical system, there are no forward orbits in $(X,G)$. 
\end{example}

\begin{definition}
Let $(X,G)$ be a CR-dynamical system. We say that  
\begin{enumerate}
\item $(X,G)$ is \emph{ \color{blue}  $1^{\oplus}$-minimal} if for each $x\in X$,  $T_G^{+}(x)\neq \emptyset$, and for each $\mathbf x\in \star_{i=1}^{\infty}G$,
$$
\Cl\Big(\Orbit(\mathbf x)\Big)=X.
$$
\item $(X,G)$ is \emph{ \color{blue}  $2^{\oplus}$-minimal }if for each $x\in X$ there is $\mathbf x\in T_G^{+}(x)$ such that 
$$
\Cl\Big(\Orbit(\mathbf x)\Big)=X.
$$
\item $(X,G)$ is \emph{ \color{blue}  $3^{\oplus}$-minimal} if for each $x\in X$,
$$
\Cl\Big(\mathcal U_G^{\oplus}(x)\Big)=X.
$$
\end{enumerate}
\end{definition}

\begin{theorem}\label{main1}
Let $(X,G)$ be a CR-dynamical system. Then the following holds.
 \begin{enumerate}
\item \label{1dva}  $(X,G)$ is $1$-minimal if and only if $(X,G)$ is $1^{\oplus}$-minimal.
\item \label{2dva} If $(X,G)$ is $1^{\oplus}$-minimal, then $(X,G)$ is $2^{\oplus}$-minimal.
\item \label{3dva} If $(X,G)$ is $2^{\oplus}$-minimal, then $(X,G)$ is $3^{\oplus}$-minimal.
\item \label{4dva} If $(X,G)$ is $3^{\oplus}$-minimal, then $(X,G)$ is $\infty$-minimal.
\end{enumerate}
\end{theorem}
\begin{proof}
Let $(X,G)$ be a $1$-minimal CR-dynamical system. To prove that $(X,G)$ is $1^{\oplus}$-minimal, let  $x\in X$.  
To prove that
$$
T_G^{+}(x)\neq \emptyset,
$$
we show first that $p_2(G)\subseteq p_1(G)$. Suppose that $p_2(G)\not \subseteq p_1(G)$ and let $x_0\in p_2(G)\setminus p_1(G)$. Then $A=\{x_0\}$  is trivially $1$-invariant in $(X,G)$---a contradiction since $A\neq X$.  Therefore, $p_2(G)\subseteq p_1(G)$.  Next, we prove that $p_2(G)=X$. Let $A=p_2(G)$ and let $x\in A$ be any point.  Since $A\subseteq p_1(G)$, it follows that  $x\in p_1(G)$. Then there is $y\in p_2(G)$ such that $(x,y)\in G$. This proves that $A$ is $1$-invariant in $(X,G)$. Since $A$ is closed in $X$ and $A\neq \emptyset$, it follows that $A=X$ since $(X,G)$ is $1$-minimal. Therefore, $p_2(G)=X$. Also, $p_1(G)=X$ follows since $p_2(G)\subseteq p_1(G)$.  Since $p_1(G)=X$,  there is a point $\mathbf x\in \star_{i=1}^{\infty}G$ such that $\pi_1(\mathbf x)=x$ and $T_G^{+}(x)\neq \emptyset$. This  completes the proof that $T_G^{+}(x)\neq \emptyset$.

Next, let $\mathbf x\in \star_{i=1}^{\infty}G$.  We show that $\Cl (\Orbit(\mathbf x))=X$. Let $A=\Cl (\Orbit(\mathbf x))$. Then $A$ is a non-empty closed subset of $X$.  Let $x\in A$ such that $x\in p_1(G)$.  We consider the following possible cases.
\begin{enumerate}
\item[(i)] $x\in \Orbit(\mathbf x)$.  Let $m$ be a positive integer such that $\pi_m(\mathbf x)=x$ and let $y=\pi_{m+1}(\mathbf x)$. Then $y\in A$ and $(x,y)\in G$.
\item[(ii)] $x\not \in \Orbit(\mathbf x)$.  Let $(z_n)$ be a sequence of points in $\Orbit(\mathbf x)$ such that $\displaystyle \lim_{n\to \infty}z_n=x$. For each positive integer $n$, let $i_n$ be a positive integer such that $\pi_{i_n}(\mathbf x)=z_n$. For each positive integer $n$,  let $y_n=\pi_{i_n+1}(\mathbf x)$, and let $(y_{j_n})$ be a convergent subsequence of the sequence $(y_n)$ and let $\displaystyle \lim_{n\to \infty}y_{j_n}=y$.  Note that for each positive integer $n$, $(z_{j_n},y_{j_n})\in G$ . Since $G$ is closed in $X\times X$ and since $\displaystyle \lim_{n\to \infty}(z_{j_n},y_{j_n})=(x,y)$, it follows that $(x,y)\in G$.  Since $A$ is closed in $X$ and since $y_{i_n}\in \Orbit(\mathbf x)$ for each positive integer $n$, it follows from $\Orbit(\mathbf x)\subseteq A$ that $y\in A$. 
\end{enumerate}
We proved that there is $y\in A$ such that $(x,y)\in G$.  It follows that $A$ is $1$-invariant in $(X,G)$ and, therefore, $A=X$. This proves that $\Cl(\Orbit(\mathbf x))=X$ and it follows that $(X,G)$ is $1^{\oplus}$-minimal. This proves the first implication of \ref{1dva}.

To prove the other implication, suppose that $(X,G)$ is $1^{\oplus}$-minimal and let $A$ be a non-empty closed subset of $X$ which is $1$-invariant in $(X,G)$.   Let $a_1\in A$ be any point.  Since $T_G^{+}(a_1)\neq \emptyset$, there is $\mathbf x_1\in T_G^{+}(a_1)$.  Choose such an element $\mathbf x_1\in T_G^{+}(a_1)$ and set $x=\pi_2(\mathbf x_1)$.  Then $(a_1,x)\in G$ and $a_1\in p_1(G)$ follows.  Since $A$ is $1$-invariant in $(X,G)$, there is a point $a_2\in A$ such that $(a_1,a_2)\in G=\star_{i=1}^{1}G$.   Fix such a point $a_2$.  Let $n>1$ be a positive integer and suppose that we have already constructed the points $a_1,a_2,a_3,\ldots ,a_n\in A$ such that $(a_1,a_2,a_3,\ldots ,a_n)\in \star_{i=1}^{n-1}G$.  Since $T_G^{+}(a_n)\neq \emptyset$, there is $\mathbf x_n\in T_G^{+}(a_n)$.  Choose such an element $\mathbf x_n\in T_G^{+}(a_n)$ and set $x=\pi_2(\mathbf x_n)$.  Then $(a_n,x)\in G$ and $a_n\in p_1(G)$ follows.  Since $A$ is $1$-invariant in $(X,G)$, there is a point $a_{n+1}\in A$ such that $(a_n,a_{n+1})\in G$.  Fix such a point $a_{n+1}$.  Let 
$$
\mathbf x=(a_1,a_2,a_3,\ldots).
$$
Then $\mathbf x\in \star_{i=1}^{\infty}G$ and $\Cl(\Orbit(\mathbf x))=X$ and, since $\Cl(\Orbit(\mathbf x))\subseteq A$, it follows that $A=X$. This proves that $(X,G)$ is $1$-minimal and we have just proved \ref{1dva}.

To prove \ref{2dva} suppose that $(X,G)$ is $1^{\oplus}$-minimal.  Let $x\in X$ be any point.  Since $(X,G)$ is $1^{\oplus}$-minimal, there is a point $\mathbf x\in \star_{i=1}^{\infty}G$ such that $\pi_1(\mathbf x)=x$.  Since $(X,G)$ is $1^{\oplus}$-minimal, it follows that $\Cl(\Orbit(\mathbf x))=X$.  Therefore,  $(X,G)$ is $2^{\oplus}$-minimal.

To prove \ref{3dva} suppose that $(X,G)$ is $2^{\oplus}$-minimal.  Let $x\in X$ be any point.  Since $(X,G)$ is a $2^{\oplus}$-minimal dynamical system, there is a point $\mathbf x_0\in \star_{i=1}^{\infty}G$ such that $\pi_1(\mathbf x_0)=x$ and $\Cl(\Orbit(\mathbf x_0))=X$. It follows from $\Orbit(\mathbf x_0)\subseteq \mathcal U_G^{\oplus}(x)$ that $\Cl(\mathcal U_G^{\oplus}(x))=X$.  Therefore,  $(X,G)$ is $3^{\oplus}$-minimal.

Finally, to prove \ref{4dva}, suppose first, that  $(X,G)$ is $3^{\oplus}$-minimal.  Let $A$ be a non-empty closed subset of $X$  such that $A$ is $\infty$-invariant in $(X,G)$.  Let $x\in A$.  Since $(X,G)$ is $3^{\oplus}$-minimal, it follows that $\Cl(\mathcal U_G^{\oplus}(x))=X$.  We show that $A=X$ by showing that $\Cl(\mathcal U_G^{\oplus}(x))\subseteq A$.  First, we show that $\mathcal U_G^{\oplus}(x)\subseteq A$.  Let $y\in \mathcal U_G^{\oplus}(x)$ and let $\mathbf x_0\in T_G^{+}(x)$ such that $y\in \Orbit(\mathbf x_0)$. Since $x\in A$ and since  $A$ is $\infty$-invariant in $(X,G)$, it follows that $y\in A$ by Observation \ref{gorane1}.  Therefore,  $\mathcal U_G^{\oplus}(x)\subseteq A$ and, since $A$ is closed in $X$, it follows that $\Cl(\mathcal U_G^{\oplus}(x))\subseteq A$. 
\end{proof}

Among other things, the following theorem says that $1$, $1^{\oplus}$, $2^{\oplus}$, $3^{\oplus}$, and $\infty$-minimality of  CR-dynamical systems is a generalization of the notion of the minimality of dynamical systems. 
\begin{theorem}
Let $(X,f)$ be a dynamical system. The following statements are equivalent.
\begin{enumerate}
\item \label{1ena} $(X,f)$ is minimal.
\item \label{2ena} $(X,\Gamma(f))$ is $1$-minimal.
\item \label{3ena} $(X,\Gamma(f))$ is $1^{\oplus}$-minimal.
\item \label{4ena} $(X,\Gamma(f))$ is $2^{\oplus}$-minimal.
\item \label{5ena} $(X,\Gamma(f))$ is $3^{\oplus}$-minimal.
\item \label{6ena} $(X,\Gamma(f))$ is $\infty$-minimal.
\end{enumerate}
\end{theorem}
\begin{proof}
Suppose that $(X,f)$ is minimal. To prove that $(X,\Gamma(f))$ is $1$-minimal, let  $A$ be a closed subset of $X$ such that $A$ is  $1$-invariant in $(X,\Gamma(f))$. By Observation \ref{obs1}, $A$ is $f$-invariant in $(X,\Gamma(f))$. Therefore,  $ A\in \{\emptyset, X\}$ since $(X,f)$ is minimal. This proves the implication from \ref{1ena} to \ref{2ena}.

The implications from \ref{2ena} to \ref{3ena}, from \ref{3ena} to \ref{4ena}, from \ref{4ena} to \ref{5ena} and from \ref{5ena} to \ref{6ena} follow from Theorem \ref{main1}.

Suppose that $(X,\Gamma(f))$ is $\infty$-minimal. To prove that $(X,f)$ is minimal, let  $A$ be a closed subset of $X$ such that $A$ is  $f$-invariant.  By Observation \ref{obs1}, $A$ is $\infty$-invariant in $(X,\Gamma(f))$. Therefore,  $ A\in \{\emptyset, X\}$ since $(X,\Gamma(f))$ is $\infty$-minimal. This proves the implication from \ref{6ena} to \ref{1ena}.
\end{proof}

\begin{theorem}\label{surjektivnost}
Let $(X,G)$ be a CR-dynamical system and let $k\in \{1,1^{\oplus},2^{\oplus},3^{\oplus},\infty\}$. If $(X,G)$ is $k$-minimal, then 
$$
p_1(G)=p_2(G)=X.
$$
\end{theorem}
\begin{proof}
Suppose that $(X,G)$ is $\infty$-minimal.  First, we show that $p_2(G)\subseteq p_1(G)$. Suppose that $p_2(G)\not \subseteq p_1(G)$ and let $x_0\in p_2(G)\setminus p_1(G)$. Then $A=\{x_0\}$ is trivially $\infty$-invariant in $(X,G)$---a contradiction.  Therefore, $p_2(G)\subseteq p_1(G)$. 

Next, we prove that $p_2(G)=X$. Let $A=p_2(G)$ and let $(x,y)\in G$ be any point such that $x\in A$. Since $A\subseteq p_1(G)$, it follows that $y\in p_2(G)$, meaning that $y\in A$. This proves that $A$ is $\infty$-invariant in $(X,G)$. Since $A$ is closed in $X$ and $A\neq \emptyset$, it follows that $A=X$ since $(X,G)$ is $\infty$-minimal. Therefore, $p_2(G)=X$. Also, $p_1(G)=X$ follows since $p_2(G)\subseteq p_1(G)$.

Next, let $k\in\{1,1^{\oplus},2^{\oplus},3^{\oplus}\}$ and suppose that $(X,G)$ is $k$-minimal. It follows from Theorem \ref{main1} that $(X,G)$ is also $\infty$-minimal. Therefore, $p_1(G)=p_2(G)=X$.
\end{proof}

In the following example,  we show that there is a $2^{\oplus}$-minimal CR-dynamical system which is not $1^{\oplus}$-minimal.
\begin{example} \label{ex2}
Let $X=[0,1]$ and let $G=([0,1]\times \{\frac{1}{2}\})\cup(\{\frac{1}{2}\}\times [0,1])$, see Figure \ref{figure1}.
To show that $(X,G)$ is not $1^{\oplus}$-minimal, let $\mathbf x=(\frac{1}{2}, \frac{1}{2}, \frac{1}{2}, \ldots)\in \star_{i=1}^{\infty}G$. Then $\Cl(\mathcal O_G^{\oplus}(\mathbf x))\neq X$. Therefore, $(X,G)$ is not $1^{\oplus}$-minimal.

To show that $(X,G)$ is $2^{\oplus}$-minimal, let $x\in X$ be any point.  We show that there is $\mathbf x\in T_G^{+}(x)$ such that $\Cl\Big(\Orbit(\mathbf x)\Big)=X$. Let $[0,1]\cap \mathbb Q=\{q_1,q_2,q_3,\ldots\}$ be the set of rationals in $[0,1]$, let $x_1=x$,  for each positive integer $n$,  let $x_{2n}=\frac{1}{2}$ and $x_{2n+1}=q_n$, and let $\mathbf x=(x_1,x_2,x_3,\ldots)$. Then $\mathbf x\in T_G^{+}(x)$ such that $\Cl\Big(\Orbit(\mathbf x)\Big)=X$.
\end{example}
In the following example,  we show that there is a $\infty$-minimal CR-dynamical system which is not $3^{\oplus}$-minimal.
\begin{example} \label{ex22}
Let $X=[0,1]$ and let $G$ be the union of the following line segments:
\begin{enumerate}
\item the line segment with endpoints $(0,\frac{1}{2})$ and $(1,1)$,
\item the line segment with endpoints $(1,0)$ and $(1,1)$,
\end{enumerate} 
see Figure \ref{figure2}.
\begin{figure}[h!]
	\centering
		\includegraphics[width=15em]{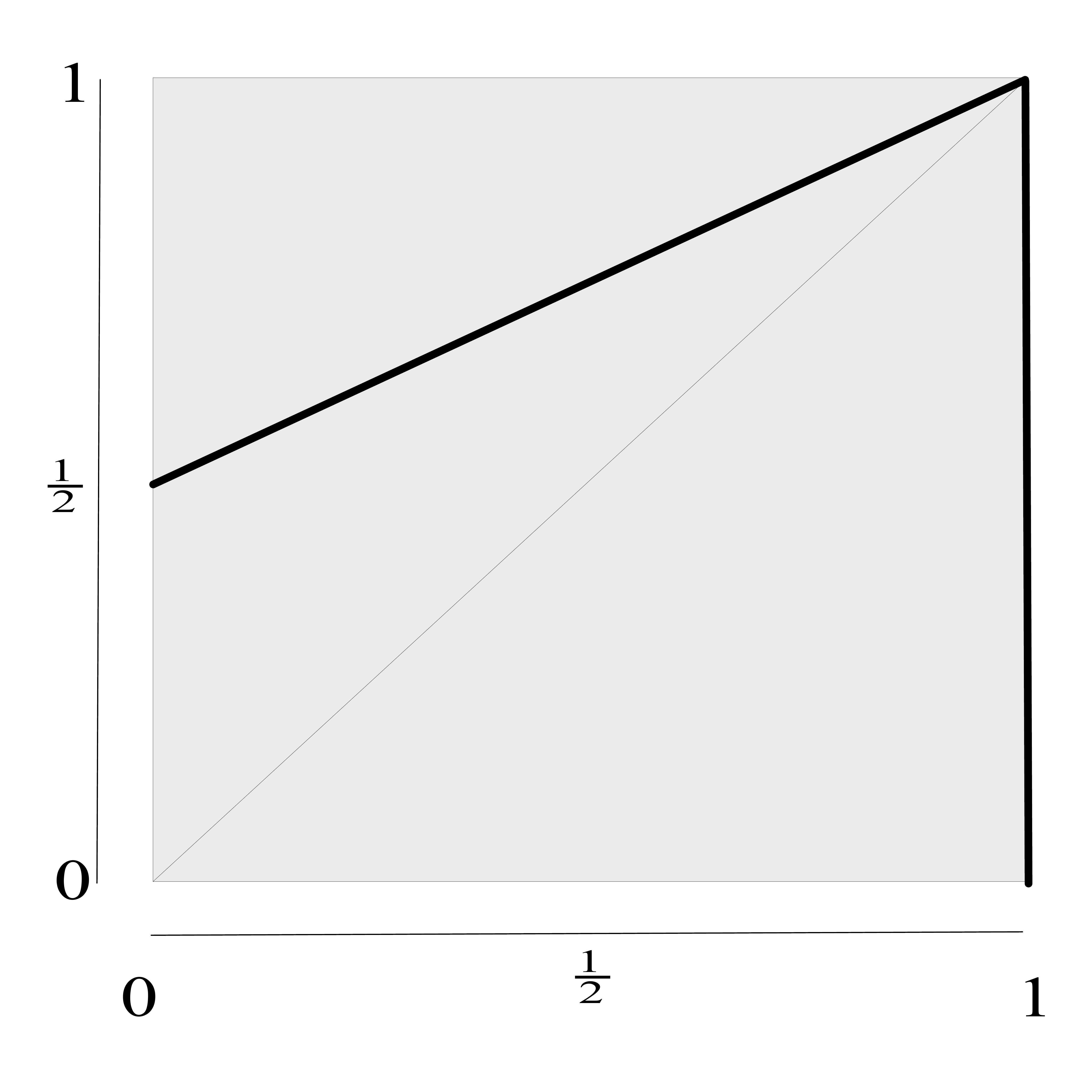}
	\caption{The relation  $G$ from Example \ref{ex22}}
	\label{figure2}
\end{figure}  
To show that $(X,G)$ is not $3^{\oplus}$-minimal, let $x_0=0$.  Then 
$$
\Cl(\mathcal U_G^{\oplus}(0))=\Cl\Big(\Big\{0,\frac{1}{2},\frac{3}{4}, \frac{7}{8},\frac{15}{16},\ldots\Big\}\Big)=\Big\{0,\frac{1}{2},\frac{3}{4}, \frac{7}{8},\frac{15}{16},\ldots\Big\}\cup\{1\}\neq X.
$$ 
Therefore, $(X,G)$ is not $3^{\oplus}$-minimal.

To show that $(X,G)$ is $\infty$-minimal, let $A$ be a non-empty closed subset of $X$ such that $A$ is $\infty$-invariant in $(X,G)$.  First, we show that $1\in A$. Since $A\neq \emptyset$, it follows that there is $x\in A$.  Choose any element $x$ in $A$.  If $x=1$, we are done.  Suppose that $x<1$ and let $f:[0,1]\rightarrow [0,1]$ be defined by $f(t)=\frac{1}{2}t+\frac{1}{2}$ for each $t\in [0,1]$.    Note that the graph of $f$ is the line segment from $(0,\frac{1}{2})$ to $(1,1)$. Since $A$ is $\infty$-invariant in $(X,G)$, it follows from Observation \ref{gorane1}  that 
$$
\orbit(x)=\{x,f(x),f^{2}(x),f^{3}(x),\ldots\}\subseteq A.
$$
Since $A$ is closed in $X$, it follows that 
$$
\Cl(\orbit(x))=\{x,f(x),f^{2}(x),f^{3}(x),\ldots\}\cup \{1\}\subseteq A.
$$
Therefore, $1\in A$.  Next, let $y\in X$ be any point. Then $(1,y)\in G$ and since $1\in A$, it follows from the fact that $A$ is $\infty$-invariant in $(X,G)$, that $y\in A$. Therefore, $A=X$.  
\end{example}
We conclude this section by stating the following problem.

{
\begin{problem}
Is there an example of a $3^{\oplus}$-minimal CR-dynamical system which is not $2^{\oplus}$-minimal?
\end{problem}}

\section{Minimality and omega limit sets}\label{s3}
Theorem \ref{main666}, where results about relations of omega limits sets in CR-dynamical systems $(X,G)$ and minimality are presented, is the main result of this section. First, we revisit omega limit sets in dynamical systems $(X,f)$.
\begin{definition}
Let $(X,f)$ be a dynamical system and let $x_0\in X$ and let $\mathbf x\in T_{\Gamma(f)}^{+}(x_0)$ be the forward orbit of $x_0$. The set 
$$
\omega_f(x_0)=\{x\in X \ | \ \textup{ there is a subsequence of  the sequence } \mathbf x  \textup{ with limit } x\}
$$   
is called \emph{ \color{blue} the omega limit set of $x_0$}. 
\end{definition}

The following is a well-known result.  We present its proof for the completeness of the paper.
\begin{theorem}\label{omomg}
Let $(X,f)$ be a dynamical system. The following statements are equivalent.
\begin{enumerate}
\item \label{tritri} $(X,f)$ is minimal.
\item \label{stiri}
For each $x\in X$, 
$$
\omega_f(x)=X.
$$
\end{enumerate}
\end{theorem}
\begin{proof}
To show that \ref{tritri} is equivalent to \ref{stiri}, suppose that for each $\mathbf x\in \star_{i=1}^{\infty}\Gamma(f)$,  $\Cl (\orbit(\mathbf x))=X$, let $x_0\in X$ be any point and let $\mathbf x_0\in T_{\Gamma(f)}^{+}(x_0)$ be the forward orbit of $x_0$. Obviously, $\omega_f(x_0)\subseteq X$.  To show that $X\subseteq \omega_f(x_0)$, let $x\in X$ and let $\mathbf x\in T_{\Gamma(f)}^{+}(x)$ be the forward orbit of $x$.  To show that $x\in \omega_f(x_0)$, we  consider the following possible cases.
\begin{enumerate}
\item There is a positive integer $n$ such that $f^n(x)=x$.  Then $x$ is a periodic point and it follows from $\Cl(\mathcal O_f^{\oplus}(\mathbf x))=X$ and from the fact that $\mathcal O_f^{\oplus}(\mathbf x)$ is finite, that $X=\mathcal O_f^{\oplus}(\mathbf x)$, so $X$ is finite.   Then $\mathcal O_f^{\oplus}(\mathbf x)=\mathcal O_f^{\oplus}(\mathbf x_0)$ and it follows that $x\in \omega_f(x_0)$.
\item For each positive integer $n$, $f^n(x)\neq x$.  First, we show that $x\in \omega_f(x_0)$. Suppose the contrary, that $x\not\in \omega_f(x_0)$. Then $x$ is not a limit point of the sequence $\mathbf x$.  Therefore, there is an open set  $U$ in $X$ such that $x\in U$ and  $U\cap \mathcal O_f^{\oplus}(\mathbf x_0)$ only contains one element.  Choose such an open set $U$ in $X$ and let $U\cap \mathcal O_f^{\oplus}(\mathbf x_0)=\{z_0\}$.  First, we show that $U\subseteq  \mathcal O_f^{\oplus}(\mathbf x_0)$.  Suppose that  $U\not \subseteq  \mathcal O_f^{\oplus}(\mathbf x_0)$. Then let $y\in U\setminus \mathcal O_f^{\oplus}(\mathbf x_0)$ and let $r_1=d(y,z_0)$, $r_2=d(y,X\setminus U)$, and let
$$
r=\min\{r_1,r_2\}.
$$
Then $r>0$ and  $B(y,\frac{r}{2})=\{z\in X \ | \ d(z,y)<\frac{r}{2}\}\subseteq X\setminus \mathcal O_f^{\oplus}(\mathbf x_0)$.  Therefore, $\mathcal O_f^{\oplus}(\mathbf x_0)$ is not dense in $X$---a contradiction.  Therefore,  $U$ is a subset of $ \mathcal O_f^{\oplus}(\mathbf x_0)$.  Since $U\cap  \mathcal O_f^{\oplus}(\mathbf x_0)$ is finite, it follows that $U$ is a finite subset of $ \mathcal O_f^{\oplus}(\mathbf x_0)$.  Let $V=\{x\}$ and let $\mathbf y\in T_{\Gamma(f)}^{+}(f(x))$ be the forward orbit of $f(x)$.  It follows that $V$ is an open set in $X$ and  that $x\not \in \mathcal O_f^{\oplus}(\mathbf y)$ since for each positive integer $n$, $f_n(x)\neq x$.  Therefore, $V\cap \mathcal O_f^{\oplus}(\mathbf y)=\emptyset$ and it follows that $\mathcal O_f^{\oplus}(\mathbf y)$ is not dense in $X$---a contradiction. Therefore,  $x\in \omega_f(x_0)$.
\end{enumerate}
We have just proved that also $X\subseteq \omega_f(x_0)$.  So, $ \omega_f(x_0)=X$ follows.

Finally, we show that \ref{stiri} implies \ref{tritri}. Suppose that for each $x\in X$, $\omega_f(x)=X$, let $x_0\in X$ be any point and let $\mathbf x_0\in T_{\Gamma(f)}^{+}(x_0)$ be the forward orbit of $x_0$.   Obviously,  $\Cl (\orbit(\mathbf x_0))\subseteq X$. To show that $X\subseteq \Cl (\orbit(\mathbf x_0))$, let $x\in X$.  Since $X=\omega_f(x_0)$, there is a subsequence $(f^{i_n}(x_0))$ of $\mathbf x_0$ such that 
$$
\lim_{n\to \infty}f^{i_n}(x_0)=x.
$$
Therefore, $x\in \Cl (\orbit(\mathbf x_0))$. This completes the proof. 
\end{proof}

Next, we generalize the notion of omega limit sets from dynamical systems to CR-dynamical systems.
\begin{definition}
Let $(X,G)$ be a CR-dynamical system and let $\mathbf x\in \star_{i=1}^{\infty}G$.  The set 
$$
\omega_G(\mathbf x)=\{x\in X \ | \ \textup{ there is a subsequence of the sequence } \mathbf x \textup{ with limit } x\}
$$
is called \emph{ \color{blue}  the  omega limit set of $\mathbf x$}.  For each  $x\in X$,  we use \emph{ \color{blue} $\psi_G(x)$} to denote the set 
$$
\psi_G(x)=\bigcup_{\mathbf x\in T_G^{+}(x)}\omega_G(\mathbf x).
$$
\end{definition}

\begin{observation}\label{gorane2}
Note that for each $\mathbf x\in \star_{i=1}^{\infty}G$, $\omega_G(\mathbf x)\subseteq \Cl(\mathcal O_G^{\oplus}(\mathbf x))$.
\end{observation}

\begin{definition}
Let $(X,G)$ be a CR-dynamical system.  We say that 
\begin{enumerate}
\item $(X,G)$ is \emph{ \color{blue} $1^{\omega}$-minimal}, if for each $x\in X$,  $T_G^{+}(x)\neq \emptyset$, and for each $\mathbf x\in \star_{i=1}^{\infty}G$,
$$
\omega_G(\mathbf x)=X.
$$
\item $(X,G)$ is \emph{ \color{blue} $2^{\omega}$-minimal}, if for each $x\in X$ there is $\mathbf x\in T_G^{+}(x)$ such that 
$$
\omega_G(\mathbf x)=X.
$$
\item $(X,G)$ is \emph{ \color{blue} $3^{\omega}$-minimal}, if for each $x\in X$,
$$
\psi_G(x)=X.
$$
\end{enumerate}
\end{definition}
\begin{observation}
Let $(X,G)$ be a CR-dynamical system. Then the following hold. 
\begin{enumerate}
\item If $(X,G)$ is $1^{\omega}$-minimal, then $(X,G)$ is $2^{\omega}$-minimal.
\item If $(X,G)$ is $2^{\omega}$-minimal, then $(X,G)$ is $3^{\omega}$-minimal.
\end{enumerate}
\end{observation}
Next, we construct an example of a CR-dynamical system such that it is $2^{\omega}$-minimal but it is not $1^{\omega}$-minimal. 
\begin{example}\label{omegaPLUS}
Let $X=[0,1]$ and let $G=([0,1]\times \{\frac{1}{2}\})\cup(\{\frac{1}{2}\}\times [0,1])$, see Figure \ref{figure1}.
Then $(X,G)$ is $2^{\omega}$-minimal but it is not $1^{\omega}$-minimal. 

To show that $(X,G)$ is not $1^{\omega}$-minimal, let $\mathbf x=(\frac{1}{2}, \frac{1}{2}, \frac{1}{2}, \ldots)\in T_G^{+}(\frac{1}{2})$. Then $\omega_G(\mathbf x)\neq X$. Therefore, $(X,G)$ is not $1^{\omega}$-minimal.

To show that $(X,G)$ is $2^{\omega}$-minimal, let $x\in X$ be any point.  We show that there is $\mathbf x\in T_G^{+}(x)$ such that $\omega_G(\mathbf x)=X$. Let $\mathbb Q$ denote the rationals and let $[0,1]\cap \mathbb Q=\{q_1,q_2,q_3,\ldots\}$, let $x_1=x$,  for each positive integer $n$,  let $x_{2n}=\frac{1}{2}$ and $x_{2n+1}=q_n$, and let $\mathbf x=(x_1,x_2,x_3,\ldots)$. Then $\mathbf x\in T_G^{+}(x)$ such that $\omega_G(\mathbf x)=X$.
\end{example}

\begin{observation}
Let $(X,f)$ be a dynamical system,  let $x\in X$ and let $\mathbf x\in T_{\Gamma(f)}^{+}(x)$ be the forward orbit of $x$.  Then 
$$
\omega_f(x)=\omega_{\Gamma(f)}(\mathbf x)=\psi_{\Gamma(f)}(x)
$$
and the following statements are equivalent.
\begin{enumerate}
\item $(X,f)$ is minimal.
\item $(X,\Gamma(f))$ is $1^{\omega}$-minimal.
\item $(X,\Gamma(f))$ is $2^{\omega}$-minimal.
\item $(X,\Gamma(f))$ is $3^{\omega}$-minimal.
\end{enumerate}
\end{observation}
\begin{theorem}\label{main666}
Let $(X,G)$ be a CR-dynamical system. Then the following hold.
\begin{enumerate}
\item \label{1tre} $(X,G)$ is $1^{\omega}$-minimal if and only if $(X,G)$ is $1^{\oplus}$-minimal. 
\item  \label{2tre}  $(X,G)$ is $2^{\omega}$-minimal if and only if $(X,G)$ is $2^{\oplus}$-minimal.
\item  \label{3tre} If $(X,G)$ is $3^{\omega}$-minimal, then $(X,G)$ is $3^{\oplus}$-minimal.
\end{enumerate}
\end{theorem}
\begin{proof}
To prove \ref{1tre},  first suppose that $(X,G)$ is $1^{\oplus}$-minimal. Clearly,  for each for each $x\in X$,  $T_G^{+}(x)\neq \emptyset$. Let $\mathbf x\in \star_{i=1}^{\infty}G$. Obviously,  $\omega_G(\mathbf x)\subseteq X$.  To prove that $X\subseteq \omega_G(\mathbf x)$, let $x\in X$. To show that $x\in \omega_G(\mathbf x)$, we treat the following possible cases.
\begin{enumerate}
\item[(i)] $x\not \in \Orbit(\mathbf x)$. Since $\Orbit(\mathbf x)$ is dense in $X$, it follows that for any open set $U$ in $X$, 
$$
U\neq \emptyset \Longrightarrow U \cap \Orbit(\mathbf x)\neq \emptyset.
$$
Therefore, for any open set $U$ in $X$, 
$$
x\in U \Longrightarrow U\cap \Orbit(\mathbf x)\neq \emptyset
$$
and since $x\not\in \Orbit(\mathbf x)$, it follows that for any open set $U$ in $X$, 
$$
x\in U \Longrightarrow (U\setminus \{x\})\cap \Orbit(\mathbf x)\neq \emptyset.
$$
Therefore, $x$ is a limit point of the sequence $\mathbf x$ and $x\in \omega_G(\mathbf x)$ follows. 
\item[(ii)] $x \in \Orbit(\mathbf x)$. Suppose that $x\not \in \omega_G(\mathbf x)$. Then there is an open set $U$ in $X$ such that $U\cap \Orbit(\mathbf x)=\{x\}$ and $\pi_k(\mathbf x)=x$ only for finitely many positive integers $k$. Let 
$$
n_0=\max\{n \ | \ n \textup{ is a positive integer such that } \pi_n(\mathbf x)=x\}
$$
and let
$$
\mathbf y=(\pi_{n_0+1}(\mathbf x),\pi_{n_0+2}(\mathbf x),\pi_{n_0+3}(\mathbf x),\ldots).
$$
Then $U\cap \Orbit(\mathbf y)=\emptyset$---a contradiction since $(X,G)$ is $1^{\oplus}$-minimal and, therefore, $\Cl(\Orbit(\mathbf y))=X$. It follows that $x\in \omega_G(\mathbf x)$.
\end{enumerate}
Next, suppose that for each $x\in X$,  $T_G^{+}(x)\neq \emptyset$, and for each $\mathbf x\in \star_{i=1}^{\infty}G$,
$$
\omega_G(\mathbf x)=X.
$$
Therefore,  for each $x\in X$,  $T_G^{+}(x)\neq \emptyset$, and by Observation \ref{gorane2},  for each $\mathbf x\in \star_{i=1}^{\infty}G$, 
$$
X=\omega_G(\mathbf x)\subseteq \Cl(\Orbit(\mathbf x))\subseteq X
$$
and $\Cl(\Orbit(\mathbf x))= X$ follows. This completes the proof of \ref{1tre}.

Next, we prove \ref{2tre}. Suppose that $(X,G)$ is $2^{\omega}$-minimal, i.e. that for each $x\in X$ there is $\mathbf x\in T_G^{+}(x)$ such that 
$$
\omega_G(\mathbf x)=X.
$$ 
To show that $(X,G)$ is $2^{\oplus}$-minimal, let $x_0\in X$ be any point and let $\mathbf x_0\in T_G^{+}(x_0)$ be such that $\omega_G(\mathbf x_0)=X$.  By Observation \ref{gorane2}, 
$$
X=\omega_G(\mathbf x_0)\subseteq  \Cl(\Orbit(\mathbf x_0))\subseteq X.
$$
 Therefore, $\Cl(\Orbit(\mathbf x_0))= X$. This completes the proof of one implication of \ref{2tre}.

 Next, suppose that $(X,G)$ is $2^{\oplus}$-minimal. To show that $(X,G)$ is $2^{\omega}$-minimal, let $x\in X$ be any point. We will construct $\mathbf x\in T_G^{+}(x)$ such that $\omega_G(\mathbf x)=X$. Let 
$$
\mathbf x_1=(x_1^1,x_2^1,x_3^1,\ldots )\in T_G^{+}(x)
$$
 such that $\Cl(\Orbit(\mathbf x_1))=X$. 
For each positive integer $n$, let $\ell_n$ be a positive integer and let $y_1^{n},y_2^{n},y_3^{n},\ldots ,y_{\ell_n}^{n}\in X$, such that 
$$
\mathcal U_n=\Big\{B(y_1^{n},\frac{1}{n}),B(y_2^{n},\frac{1}{n}),B(y_2^{n},\frac{1}{n}),\ldots ,B(y_{\ell_n}^{n},\frac{1}{n})\Big\}
$$
is a finite open cover for $X$.  
We follow the following steps.

{\bf Step 1.} Let $m_1$ be a positive integer such that for each $i\in \{1,2,3,\ldots ,\ell_1\}$,
$$
\{x_1^{1},x_2^{1},x_3^{1},\ldots ,x_{m_1}^{1}\}\cap B(y_i^{1},1)\neq \emptyset.
$$

{\bf Step 2.} Let 
$$
\mathbf x_2=(x_1^{2},x_2^{2},x_3^{2},\ldots )\in T_G^{+}(x_{m_1}^{1})
$$
 and let $m_2$ be a positive integer such that for each $i\in \{1,2,3,\ldots ,\ell_{2}\}$,
$$
\{x_1^{2},x_2^{2},x_3^{2},\ldots ,x_{m_2}^{2}\}\cap B(y_i^{2},\frac{1}{2})\neq \emptyset.
$$

{\bf Step 3.} Let 
$$
\mathbf x_3=(x_1^{3},x_2^{3},x_3^{3},\ldots )\in T_G^{+}(x_{m_2}^{2})
$$
 and let $m_3$ be a positive integer such that for each $i\in \{1,2,3,\ldots ,\ell_{3}\}$,
$$
\{x_1^{3},x_2^{3},x_3^{3},\ldots ,x_{m_3}^{3}\}\cap B(y_i^{3},\frac{1}{3})\neq \emptyset.
$$

We continue  inductively. For each positive integer $j$, the step $j$ is as follows.

{\bf Step $\mathbf j$.} Let 
$$
\mathbf x_j=(x_1^{j},x_2^{j},x_3^{j},\ldots )\in T_G^{+}(x_{m_{j-1}}^{j-1})
$$
 and let $m_j$ be a positive integer such that for each $i\in \{1,2,3,\ldots ,\ell_{j}\}$,
$$
\{x_1^{j},x_2^{j},x_3^{j},\ldots ,x_{m_j}^{j}\}\cap B(y_i^{j},\frac{1}{j})\neq \emptyset.
$$

Finally, let 
$$
\mathbf x=(x_1^{1},x_2^{1},x_3^{1},\ldots ,x_{m_1}^{1}=x_1^{2},x_2^{2},x_3^{2},\ldots ,x_{m_2}^{2}=x_1^{3},x_2^{3},x_3^{3},\ldots ,x_{m_3}^{3},\ldots).
$$
Then $\mathbf x\in T_G^{+}(x)$ such that $\omega_G(\mathbf x)=X$.

Finally, we prove \ref{3tre}.  Suppose that for each $x\in X$, $ \psi_G(x)=X$. To prove that $(X,G)$ is $3^{\oplus}$-minimal, let $x_0\in X$ be any point and we show that $\Cl(\Orbitt(x_0))=X$.  Obviously, $\Cl(\Orbitt(x_0))\subseteq X$. To show that $X\subseteq  \Cl(\Orbitt(x_0))$, let $x\in X$.  Then $x\in \psi_G(x_0)$.  Let $\mathbf x_0\in T_G^{+}(x_0)$ such that $x\in \omega_G(\mathbf x_0)$.  Since $\omega_G(\mathbf x_0)\subseteq \Cl(\Orbit(\mathbf x_0))$, it follows that  $x\in \Cl(\Orbit(\mathbf x_0))$.  Since $\Orbit(\mathbf x_0)\subseteq \Orbitt(x_0)$, it follows that $\Cl(\Orbit(\mathbf x_0))\subseteq \Cl(\Orbitt(x_0))$ and, therefore,
$$
x\in \Cl(\Orbitt(x_0)).
$$
\end{proof}
We conclude this section by stating the following problems.

{
\begin{problem}
Is there an example of a $3^{\omega}$-minimal CR-dynamical system which is not $2^{\omega}$-minimal?
\end{problem}}

{
\begin{problem} Is there an example of a $3^{\oplus}$-minimal CR-dynamical system which is not $3^{\omega}$-minimal?
\end{problem}}

\section{Backward minimal dynamical systems with closed relations}\label{s4}
In this section, we define backward dynamical systems with closed relations. 
\begin{definition}
Let $(X,G)$ be a CR-dynamical system and let $A\subseteq X$. We say that the set $A$ is 
\begin{enumerate}
\item \emph{ \color{blue}  $1$-backward invariant in $(X,G)$}, if for each $y\in A$,
$$
y\in p_2(G) \Longrightarrow \textup{ there is  } x\in A \textup{ such that } (x,y)\in G.
$$
\item \emph{ \color{blue}  $\infty$-backward invariant in $(X,G)$}, if for each $(x,y)\in G$,
$$
y\in A \Longrightarrow x\in A.
$$
\end{enumerate} 
\end{definition}
\begin{observation}\label{ema}
Let $(X,G)$ be a CR-dynamical system and let $A\subseteq X$.  Note that
\begin{enumerate}
\item $A$ is $1$-backward invariant in $(X,G)$ if and only if $A$ is $1$-invariant in $(X,G^{-1})$.
\item$A$ is $\infty$-backward invariant in $(X,G)$ if and only if $A$ is $\infty$-invariant in $(X,G^{-1})$.
\end{enumerate} 
\end{observation}
\begin{proposition}\label{mimika11}
Let $(X,G)$ be a CR-dynamical system and let $A\subseteq X$. If $A$ is $\infty$-backward invariant in $(X,G)$, then $A$ is $1$-backward invariant in $(X,G)$.
\end{proposition}
\begin{proof}
The proposition follows from Proposition \ref{mimika} and Observation \ref{ema}.
\end{proof}
\begin{example} \label{ex1b}
Let $X=[0,1]$ and let $G=([0,1]\times \{\frac{1}{2}\})\cup(\{\frac{1}{2}\}\times [0,1])$, see Figure \ref{figure1}.
The set $A=\{\frac{1}{2}\}$ is $1$-backward invariant in $(X,G)$ but it is not $\infty$-backward invariant in $(X,G)$.
\end{example}

\begin{definition}
Let $(X,G)$ be a CR-dynamical system. We say that  
\begin{enumerate}
\item $(X,G)$ is  \emph{ \color{blue}  $1$-backward minimal} if for each closed subset $A$ of $X$,
$$
A \textup{ is } {1}\textup{-backward invariant in } (X,G) \Longrightarrow A\in \{\emptyset, X\}.
$$
\item $(X,G)$ is \emph{ \color{blue}  $\infty$-backward minimal} if for each closed subset $A$ of $X$,
$$
A \textup{ is } \infty\textup{-backward invariant in } (X,G) \Longrightarrow A\in \{\emptyset, X\}.
$$
\end{enumerate}
\end{definition}
\begin{observation}\label{333bbb}
Let $(X,G)$ be a CR-dynamical system and let $k\in \{1,\infty\}$. Then the following holds. 
$$
(X,G) \textup{ is } k\textup{-backward minimal } \Longleftrightarrow  (X,G^{-1}) \textup{ is } k\textup{-minimal}.
$$
\end{observation}

\begin{theorem}\label{judyKkkk}
Let $(X,G)$ be a CR-dynamical system.  If $(X,G)$ is $1$-backward minimal,  then $(X,G)$ is $\infty$-backward minimal. 
\end{theorem}
\begin{proof}
Let $(X,G)$ be $1$-backward minimal and let $A$ be a non-empty closed subset of $X$ such that $A$ is $\infty$-backward invariant in $(X,G)$.  Then $A$ is $1$-backward invariant in $(X,G)$ by Proposition \ref{mimika11}. Therefore, $A=X$ and it follows that $(X,G)$ is $\infty$-backward minimal.
\end{proof}
Note that Example \ref{goranH} is an example of a $\infty$-backward minimal CR-dynamical system which is not $1$-backward minimal. In Theorem \ref{666}, we show (using backward orbits that are defined in Section \ref{s5}) that for a CR-dynamical system $(X,G)$, the following holds:
$$
(X,G) \text{ is }1\text{-backward minimal} \Longleftrightarrow (X,G) \text{ is }1\text{-minimal.}
$$ 
The following example gives a CR-dynamical system which is $\infty$-minimal but is not $\infty$-backward minimal. 
\begin{example}\label{bluy}
Let $X=[0,1]$ and let $G$ be the union of the following line segments:
\begin{enumerate}
\item the line segment with endpoints $(0,\frac{1}{2})$ and $(1,1)$,
\item the line segment with endpoints $(1,0)$ and $(1,1)$,
\end{enumerate} 
see Figure \ref{figure2}. We proved that $(X,G)$ is $\infty$-minimal in Example \ref{ex22}.

To show that $(X,G)$ is not $\infty$-backward minimal, let 
$$
A=\{0,\frac{1}{2}, \frac{3}{4}, \frac{7}{8},\ldots \}\cup\{1\}.
$$
Then $A$ is a non-empty closed subset of $X$. Let $(x,y)\in G$ such that $y\in A$. If $y=0$, then $x=1$ and, therefore $x\in A$. If $y=\frac{1}{2}$, then $x=0$ and, therefore $x\in A$. If $y=1$, then $x=1$ and, therefore $x\in A$. If $y=\frac{2^{n+1}-1}{2^{n+1}}$ for some positive integer $n$, then $x=\frac{2^{n}-1}{2^{n}}$ or $x=1$, therefore, $x\in A$. This proves that $A$ is $\infty$-backward invariant. Since $A\neq X$, it follows that $(X,G)$ is not $\infty$-backward minimal.
\end{example}
The following example gives a CR-dynamical system which is $\infty$-backward minimal but is not $\infty$-minimal. 
\begin{example}
Let $X=[0,1]$ and let $H$ be the union of the following line segments:
\begin{enumerate}
\item the line segment with endpoints $(0,\frac{1}{2})$ and $(1,1)$,
\item the line segment with endpoints $(1,0)$ and $(1,1)$,
\end{enumerate} 
 and let $G=H^{-1}$. Then, using a similar approach as in Example \ref{bluy}, one can easily prove that $(X,G) $ is $\infty$-backward minimal but is not $\infty$-minimal. 
\end{example}

\section{Minimality and backward orbits}\label{s5}
First, we visit the dynamical systems and revisit a well-known result saying that a dynamical system $(X,f)$ is minimal if and only if $f$ is surjective and every backward orbit in $(X,f)$ is dense in $X$.
\begin{definition}
Let $(X,f)$ be a dynamical system and let $x_0\in X$. We use  \emph{ \color{blue} $T_f^{-}(x_0)$} to denote the set 
$$
T_f^{-}(x_0)=\{\mathbf x\in \star_{i=1}^{\infty}\Gamma(f)^{-1} \  |  \  \pi_1(\mathbf x)=x_0\}.
$$
\end{definition}
\begin{definition}
Let $(X,f)$ be a dynamical system, let $x_0\in X$ be any point  and let 
$$
\mathbf x=(x_1,x_2,x_3,\ldots) \in \star_{i=1}^{\infty}\Gamma(f)^{-1}.
$$
The sequence $\mathbf x$ is called \emph{ \color{blue} a backward orbit of $x_0$}, if  $\pi_1(\mathbf x)=x_0$.  We use  \emph{ \color{blue} $\orbitminus(\mathbf x)$ } to denote the set  
$$
\orbitminus(\mathbf x)=\{x_1,x_2,x_3,\ldots\}
$$
\end{definition}

The following is a well-known result, see \cite{KST,M} for more details.  We present its proof for the completeness of the paper.
\begin{theorem}\label{minimalback}
Let $(X,f)$ be a dynamical system. The following statements are equivalent.
\begin{enumerate}
\item \label{tritritri} $(X,f)$ is minimal.
\item \label{pet} For each $x\in X$,  $T_f^{-}(x)\neq \emptyset$ and for each $\mathbf x\in \star_{i=1}^{\infty}\Gamma(f)^{-1}$, 
$$
\Cl(\orbitminus(\mathbf x))=X.
$$
\end{enumerate}
\end{theorem}
\begin{proof}
Let $(X,f)$ be a minimal dynamical system.  Then $f$ is surjective and, therefore,  for each $x\in X$,  $T_f^{-}(x)\neq \emptyset$.  Let $\mathbf x\in \star_{i=1}^{\infty}\Gamma(f)^{-1}$.  We show that
$$
\Cl(\orbitminus(\mathbf x))=X.
$$
Let $A$ be the set of all limit points of the sequence $\mathbf x$.  Then $A\neq \emptyset$ and 
$$
A\subseteq \Cl(\orbitminus(\mathbf x)).
$$ 
To show that $A$ is closed in $X$, let $(s_n)$ be a sequence in $A$ and let $s\in X$ such that $\displaystyle \lim_{n\to \infty}s_n=s$. To show that $s\in A$, let $U$ be an open set in $X$ such that $s\in U$. Let $n_0$ be a positive integer such that $s_{n_0}\in U$. Then infinitely many terms of the sequence $\mathbf x$ are in $U$, since $s_{n_0}$ is a limit point of the sequence.  Therefore, $s\in A$ and this proves that $A$ is closed in $X$.  Next, we show that $f(A)\subseteq A$.  Let $s\in A$ and let $(x_{i_n})$ be a subsequence of the sequence $\mathbf x$ such that $\displaystyle \lim_{n\to \infty} x_{i_n}=s$. Then, since $f$ is continuous,
$$
f(s)=\lim_{n\to\infty}f(x_{i_n})=\lim_{n\to\infty}x_{i_n-1}\in A.
$$    
It follows that $f(A)\subseteq A$.  Since $(X,f)$ is minimal, it follows that $A=X$. Therefore,  $\Cl(\orbitminus(\mathbf x))=X$.

Next, suppose that for each $x\in X$,  $T_f^{-}(x)\neq \emptyset$ and for each $\mathbf x\in \star_{i=1}^{\infty}\Gamma(f)^{-1}$, 
$$
\Cl(\orbitminus(\mathbf x))=X.
$$
To show that $(X,f)$ is minimal, let $A$ be a non-empty closed subset of $X$ such that $f(A)=A$.  Suppose that $A\neq X$.  Then $\star_{i=1}^{\infty}\Gamma(f|_A)^{-1}\neq \emptyset$. Let $\mathbf x\in \star_{i=1}^{\infty}\Gamma(f|_A)^{-1}$. Then 
$$
\Cl(\orbitminus(\mathbf x))\subseteq A\neq X,
$$
which is a contradiction. Therefore, $(X,f)$ is minimal.
\end{proof}

Next, we generalize the notion of backward orbits in $(X,f)$ to the notion of backward orbits in $(X,G)$.

\begin{definition}
Let $(X,G)$ be a CR-dynamical system and let $x_0\in X$. We use \emph{ \color{blue}  $T_G^{-}(x_0)$} to denote the set 
$$
T_G^{-}(x_0)=\{\mathbf x\in \star_{i=1}^{\infty}G^{-1} \ | \  \pi_1(\mathbf x)=x_0\}.
$$
\end{definition}
\begin{definition}
Let $(X,G)$ be a CR-dynamical system, let $\mathbf x\in \star_{i=1}^{\infty}G^{-1}$, and let $x_0\in X$. We 
 \begin{enumerate}
 \item say that $\mathbf x$ is \emph{ \color{blue}  a backward  orbit of $x_0$ in $(X,G)$}, if $\pi_1(\mathbf x)=x_0$. 
 \item use \emph{ \color{blue}  $\mathcal O_G^{\ominus}(\mathbf x)$} to denote the set
 $$
 \mathcal O_G^{\ominus}(\mathbf x)=\{\pi_k(\mathbf x) \ | \  k \textup{ is a positive integer}\}.
 $$
 \item use \emph{ \color{blue}  $\mathcal U_G^{\ominus}(x_0)$} to denote the set
 $$
\mathcal U_G^{\ominus}(x_0)=\bigcup_{\mathbf x\in T_G^{-}(x_0)}\mathcal O_G^{\ominus}(\mathbf x).
 $$
 \end{enumerate}
\end{definition}
\begin{definition}
Let $(X,G)$ be a CR-dynamical system. We say that  
\begin{enumerate}
\item $(X,G)$ is \emph{ \color{blue}  $1^{\ominus}$-minimal} if for each $x\in X$,  $T_G^{-}(x)\neq \emptyset$, and for each $\mathbf x\in \star_{i=1}^{\infty}G^{-1}$,
$$
\Cl\Big(\Orbitminus(\mathbf x)\Big)=X.
$$
\item $(X,G)$ is \emph{ \color{blue}  $2^{\ominus}$-minimal }if for each $x\in X$ there is $\mathbf x\in T_G^{-}(x)$ such that 
$$
\Cl\Big(\Orbitminus(\mathbf x)\Big)=X.
$$
\item $(X,G)$ is \emph{ \color{blue}  $3^{\ominus}$-minimal} if for each $x\in X$,
$$
\Cl\Big(\mathcal U_G^{\ominus}(x)\Big)=X.
$$
\end{enumerate}
\end{definition}
\begin{observation}\label{333}
Let $(X,G)$ be a CR-dynamical system and let $k\in \{1,2,3\}$. Then the following holds. 
$$
(X,G) \textup{ is } k^{\ominus}\textup{-minimal } \Longleftrightarrow  (X,G^{-1}) \textup{ is } k^{\oplus}\textup{-minimal}.
$$
\end{observation}
\begin{theorem}\label{main1b}
Let $(X,G)$ be a CR-dynamical system. Then the following hold.
 \begin{enumerate}
\item \label{1dvab}  $(X,G)$ is $1$-backward minimal if and only if $(X,G)$ is $1^{\ominus}$-minimal.
\item \label{2dvab} If $(X,G)$ is $1^{\ominus}$-minimal, then $(X,G)$ is $2^{\ominus}$-minimal.
\item \label{3dvab} If $(X,G)$ is $2^{\ominus}$-minimal, then $(X,G)$ is $3^{\ominus}$-minimal.
\item \label{4dvab} If $(X,G)$ is $3^{\ominus}$-minimal, then $(X,G)$ is $\infty$-backward minimal.
\end{enumerate}
\end{theorem}
\begin{proof}
The proof is analogous to the proof of Theorem \ref{main1}. We leave the details to the reader.
\end{proof}
Using Observations \ref{333} and \ref{333bbb}, one can easily conclude that Example \ref{ex2} is also an example of a $2^{\ominus}$-minimal CR-dynamical system, which is not $1^{\ominus}$-minimal and that Example \ref{ex22} is also an example of a $\infty$-backward minimal CR-dynamical system, which is not $3^{\ominus}$-minimal.

\begin{theorem}\label{surjektivnostb}
Let $(X,G)$ be a CR-dynamical system. If $(X,G)$ is $1$-backward minimal, $\infty$-backward minimal or $k^{\ominus}$-minimal for some  $k\in \{1,2,3\}$, then  
$$
p_1(G)=p_2(G)=X.
$$
\end{theorem}
\begin{proof}
The theorem follows from Theorem \ref{surjektivnost} and Observations \ref{333bbb} and \ref{333}.
\end{proof}

\begin{theorem}\label{666}
Let $(X,G)$ be a CR-dynamical system.  The following statements are equivalent.
\begin{enumerate}
\item \label{GIJR1} $(X,G)$ is $1^{\ominus}$-minimal if and only if $(X,G)$ is $1^{\oplus}$-minimal. 
\item \label{GIJR4} $(X,G)$ is $1$-backward minimal if and only if $(X,G)$ is $1$-minimal. 
\end{enumerate}
\end{theorem}
\begin{proof}
First, we prove \ref{GIJR1}. Suppose that $(X,G)$ is $1^{\oplus}$-minimal.  Then $p_2(G)=X$ and, therefore,  for each $x\in X$,  $T_G^{-}(x)\neq \emptyset$.  Let $\mathbf x\in \star_{i=1}^{\infty}G^{-1}$.  Then we show that
$$
\Cl(\Orbitminus(\mathbf x))=X.
$$
Let $A$ be the set of all limit points of the sequence $\mathbf x$.  Then $A\neq \emptyset$, $A$ is closed in $X$, and 
$$
A\subseteq \Cl(\Orbitminus(\mathbf x)).
$$ 
We show that $A$ is $1$-invariant in $(X,G)$. Let $x\in A\cap p_1(G)=A$ and let $(x_{i_n})$ be a subsequence of the sequence $\mathbf x$ such that $\displaystyle \lim_{n\to \infty} x_{i_n}=x$.  Let $(s,t)$ be any limit point of the sequence $(x_{i_n},x_{i_n-1})$. Then $s=x$ and, let $y=t$. Since $G$ is closed in $X\times X$, $(x,y)\in G$ and, since $A$ is closed, it follows that $y\in A$.  We have just proved that $A$ is $1$-invariant in $(X,G)$.  Since $(X,G)$ is $1^{\oplus}$-minimal, it is also $1$-minimal by Theorem \ref{main1},  and it follows that $A=X$. Therefore,  $\Cl(\Orbitminus(\mathbf x))=X$.

Next, suppose that $(X,G)$ is $1^{\ominus}$-minimal. 
To show that $(X,G)$ is $1^{\oplus}$-minimal,  let $x\in X$ and let  $\mathbf x\in \star_{i=1}^{\infty}G$. We show that $T_G^{+}(x)\neq \emptyset$ and that $\Cl\Big(\Orbit(\mathbf x)\Big)=X$.  By Theorem \ref{surjektivnostb}, $p_1(G)=X$ and $T_G^{+}(x)\neq \emptyset$ follows.  To show that $\Cl\Big(\Orbit(\mathbf x)\Big)=X$, let $A$ be the set of all limit points of the sequence $\mathbf x$.  Then $A\neq \emptyset$, $A$ is closed in $X$, and 
$$
A\subseteq \Cl(\Orbit(\mathbf x)).
$$ 
We show that $A$ is $1$-backward invariant in $(X,G)$. Let $y\in A\cap p_2(G)=A$ and let $(x_{i_n})$ be a subsequence of the sequence $\mathbf x$ such that $\displaystyle \lim_{n\to \infty} x_{i_n}=y$.  Let $(s,t)$ be any limit point of the sequence $(x_{i_n+1},x_{i_n})$. Then $y=t$ and, let $x=s$. Since $G$ is closed in $X\times X$, $(x,y)\in G$ and, since $x$ is a limit point of $\mathbf x$, it follows that $x\in A$.  We have just proved that $A$ is $1$-backward invariant in $(X,G)$.  Since $(X,G)$ is $1^{\ominus}$-minimal, it is also $1$-backward minimal by Theorem \ref{main1b},  and it follows that $A=X$. Therefore,  $\Cl(\Orbit(\mathbf x))=X$. This completes the proof of \ref{GIJR1}. Note that this also proves \ref{GIJR4} since $(X,G)$ is $1$-minimal if and only if $(X,G)$ is $1^{\oplus}$-minimal by Theorem \ref{main1}, and since $(X,G)$ is $1$-backward minimal if and only if $(X,G)$ is $1^{\ominus}$-minimal by Theorem \ref{main1b}. 
\end{proof}

\begin{observation}
Note that in Theorem \ref{main1}, we have proved that $(X,G)$ is $1^{\oplus}$-minimal if and only if $(X,G)$ is $1$-minimal. It follows from Theorem \ref{666} that the following statements are equivalent.
\begin{enumerate}
\item $(X,G)$ is $1^{\oplus}$-minimal. 
\item  $(X,G)$ is $1^{\ominus}$-minimal. 
\item $(X,G)$ is $1$-minimal. 
\item  $(X,G)$ is $1$-backward minimal.
\end{enumerate}
\end{observation}

 Note that so far, we have not presented an example of a closed relation $G$ on $[0,1]$ such that $([0,1],G)$ is $1$-minimal. Also, note that all the closed relations $G$ on $[0,1]$ that are presented in our examples, contain a vertical or a horizontal line.  Example \ref{tistile} is an example of a closed relation $G$ on $[0,1]$ such that $([0,1],G)$ is $1$-minimal and $G$ does not contain a vertical or a horizontal line. We use Theorem \ref{tatale} in its construction.
\begin{theorem}\label{tatale}
Let $(X,G)$ be a CR-relation such that $p_1(G)=p_2(G)=X$ and let $\sigma_G:\star_{i=1}^{\infty}G^{-1}\rightarrow \star_{i=1}^{\infty}G^{-1}$ be the shift map 
$$
\sigma_G(x_1,x_2,x_3,\ldots )=(x_2,x_3,\ldots)
$$
for each $(x_1,x_2,x_3,\ldots)$. If $(\star_{i=1}^{\infty}G^{-1},\sigma_G)$ is minimal, then $(X,G)$ is $1$-minimal.
\end{theorem}
\begin{proof}
We show that $(X,G)$ is $1$-backward minimal. Let $A$ be a non-empty closed subset of $X$ such that $A$ is $1$-backward invariant. Also, let 
$$
B=\Big(\prod_{i=1}^{\infty}A\Big)\cap \Big(\star_{i=1}^{\infty}G^{-1}\Big).
$$
Since $A$ is $1$-backward invariant, $B$ is non-empty. Note, that $B$ is also a closed subset of $\star_{i=1}^{\infty}G^{-1}$ such that $\sigma_G(B)\subseteq B$. Since $(\star_{i=1}^{\infty}G^{-1},\sigma_G)$ is minimal, it follows that $B=\star_{i=1}^{\infty}G^{-1}$. Therefore, 
$$
\star_{i=1}^{\infty}G^{-1}\subseteq \prod_{i=1}^{\infty}A.
$$
Since $p_1(G)=p_2(G)=X$, it follows that 
$$
X=\pi_1(\star_{i=1}^{\infty}G^{-1})=\pi_1(B)\subseteq \pi_1(\prod_{i=1}^{\infty}A)=A.
$$
Therefore, $(X,G)$ is $1$-backward minimal. 
By Theorem \ref{666}, $(X,G)$ is $1$-minimal.
\end{proof}

\begin{example}\label{tistile}
Let $\lambda$ be an irrational number in $(0,1)$ and let $G$ be the union of the following line segments in $[0,1]\times [0,1]$:
\begin{figure}[h!]
	\centering
		\includegraphics[width=15em]{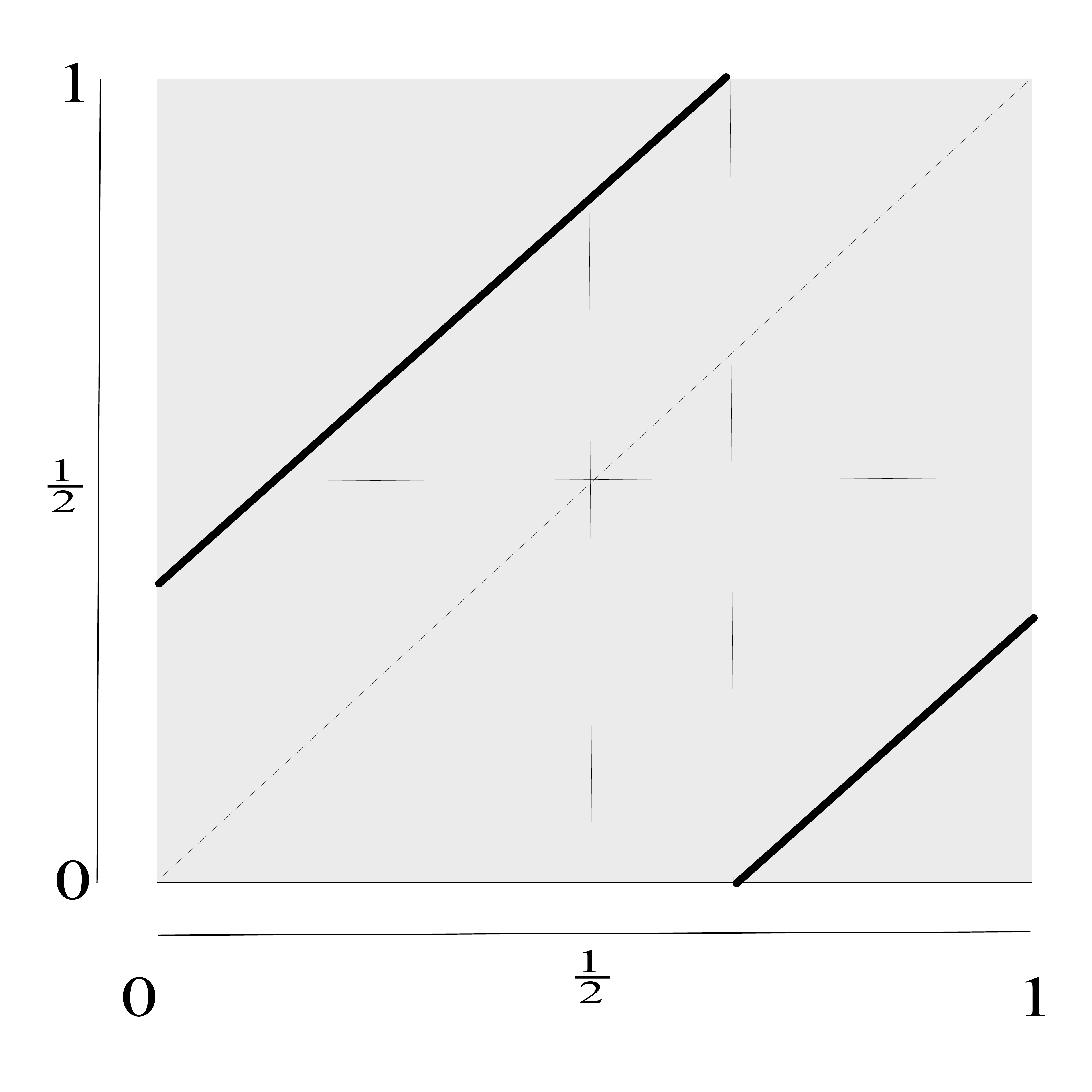}
	\caption{The relation  $G$ from Example \ref{tistile}}
	\label{figgure}
\end{figure} 
\begin{enumerate}
\item the line segment from $(0,\lambda)$ to $(1-\lambda,1)$ and
\item the line segment from $(1-\lambda,0)$ to $(1,\lambda)$,
\end{enumerate}
see Figure \ref{figgure}. 
Then $(\star_{i=1}^{\infty}G^{-1},\sigma_G)$ is minimal; this follows from the proof of \cite[Theorem 3.4, page 103]{KK}. By Theorem \ref{tatale}, $([0,1],G)$ is $1$-minimal.
\end{example}

In the following example, we demonstrate that there is a $2^{\oplus}$-minimal CR-dynamical system $(X,G)$ which is not $2^{\ominus}$-minimal.
\begin{example}\label{Rene2}
Let $X=[0,1]$ and and let $G=A\cup B\cup C$, where $A$ is a line segment from $(0,\frac{1}{2})$ to $(1,\frac{1}{2})$, $B$ is the line segment from $(0,0)$ to $(1,1)$, and $C$ is defined as follows. 

Let $d_1=\frac{1}{2}$, let $d_{10}=\frac{1}{2^{2}}$ and $d_{11}=\frac{3}{2^{2}}$, and let $d_{100}=\frac{1}{2^{3}}$, $d_{101}=\frac{3}{2^{3}}$, $d_{110}=\frac{5}{2^{3}}$ and $d_{111}=\frac{7}{2^{3}}$.  Let $n$ be a positive integer and suppose that for any  
$$
\mathbf s \in \{s_1s_2s_3\ldots s_n \ | \ s_1=1, s_2,s_3,s_4,\ldots ,s_n\in \{0,1\}\},
$$
we have already defined  $d_{\mathbf s}$ to be $d_{\mathbf s}=\frac{k}{2^n}$ for some $k\in \{1,3,5,7,\ldots ,2^n-1\}$. Then we define  $d_{\mathbf s0}$ and $d_{\mathbf s1}$ as follows.  If $k=1$ then $d_{\mathbf s0}=\frac{1}{2^{n+1}}$ and $d_{\mathbf s1}=\frac{3}{2^{n+1}}$, if $k=3$ then $d_{\mathbf s0}=\frac{5}{2^{n+1}}$ and $d_{\mathbf s1}=\frac{7}{2^{n+1}}$,  $\ldots$, and if $k=2^n-1$ then $d_{\mathbf s0}=\frac{2^{n+1}-3}{2^{n+1}}$ and $d_{\mathbf s1}=\frac{2^{n+1}-1}{2^{n+1}}$.

For each positive integer $n$, let $\mathcal S_n=\{s_1s_2s_3\ldots s_n \ | \ s_1=1, s_2,s_3,s_4,\ldots ,s_n\in \{0,1\}\}$ and let $\mathcal S=\bigcup_{n=1}^{\infty}\mathcal S_n$.
Then we define the set $C$ as 
$$
C=\bigcup_{\mathbf s\in \mathcal S} \Big(\{d_{\mathbf s}\}\times\{d_{\mathbf s0},d_{\mathbf s1}\}\Big),
$$
see Figure \ref{Rene}, 
\begin{figure}[h!]
	\centering
		\includegraphics[width=15em]{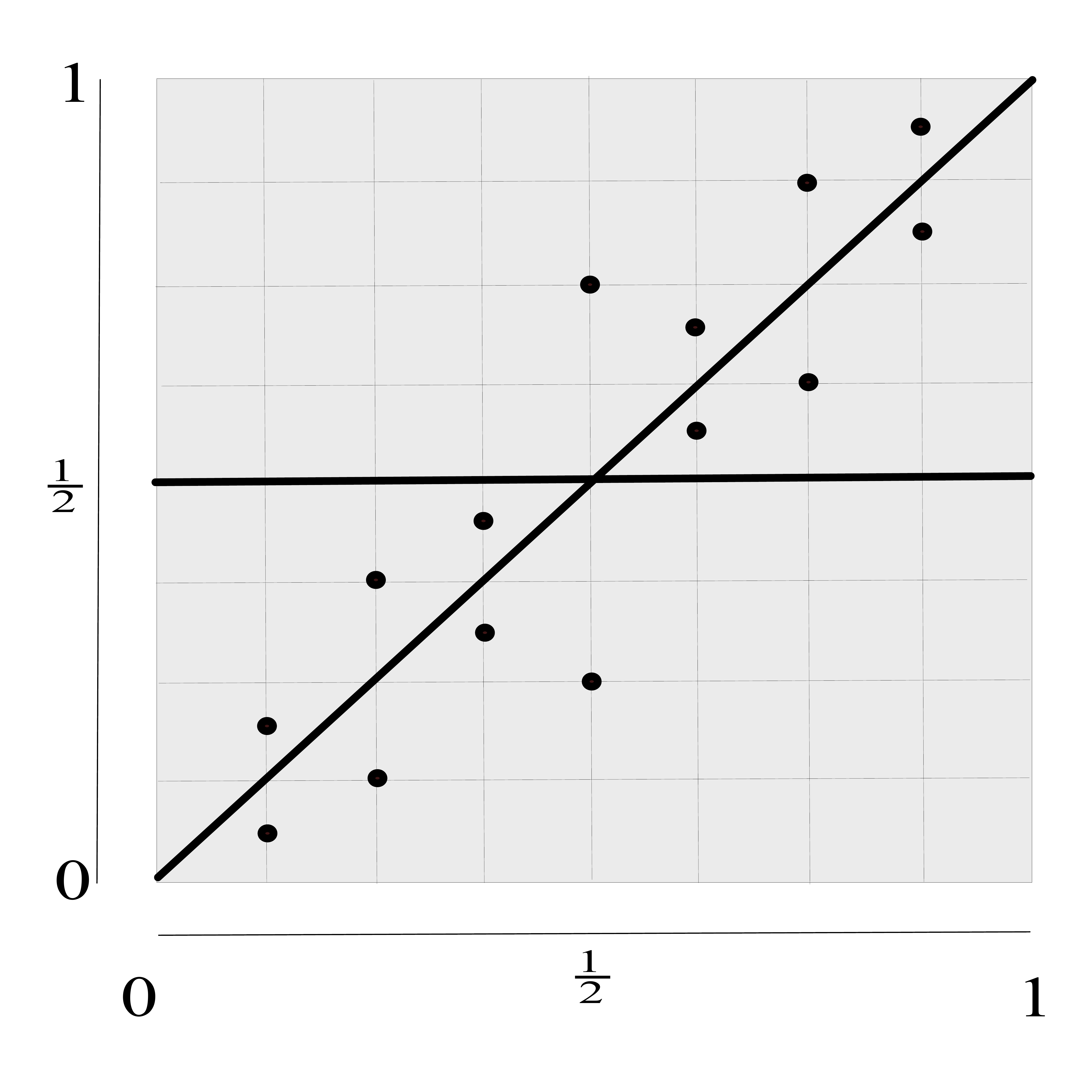}
	\caption{The construction of the set $C$ }
	\label{Rene}
\end{figure} 
where the construction of the set $C$ is presented -- in particular,  together with the sets $A$ and $B$, the set $\displaystyle \bigcup_{\mathbf s\in \mathcal S_1\cup  \mathcal S_2\cup  \mathcal S_3} \Big(\{d_{\mathbf s}\}\times\{d_{\mathbf s0},d_{\mathbf s1}\}\Big)$ is also pictured in the figure.
Then $(X,G)$ is $2^{\oplus}$-minimal (since for any $x\in [0,1]$, $(x,\frac{1}{2})\in G$ and therefore,  there is $\mathbf x\in T_G^{+}(x)$ such that $\Cl(\Orbit(\mathbf x))=X$) but it is not $2^{\ominus}$-minimal (note that $T_G^{-}(1)=\{(1,1,1,1,\ldots)\}$ and, therefore, for any $\mathbf x\in T_G^{-}(1)$, $\Cl(\mathcal O_G^{\ominus}(\mathbf x))\neq X$).
\end{example}
Note that 
\begin{enumerate}
\item $(X,G)$ from Example \ref{Rene2} is also an example of a $3^{\oplus}$-minimal  CR-dynamical system which is not $3^{\ominus}$-minimal, and
\item if $(X,G)$ is the CR-dynamical system from Example \ref{Rene2}, then $(X,G^{-1})$ is an example of a $2^{\ominus}$-minimal ($3^{\ominus}$-minimal) CR-dynamical system which is not $2^{\oplus}$-minimal ($3^{\oplus}$-minimal).
\end{enumerate}
We conclude the section by stating the following open problem.
{
\begin{problem}
 Is there an example of a $3^{\ominus}$-minimal CR-dynamical system which is not $2^{\ominus}$-minimal?
\end{problem}}
\section{Minimality and alpha limit sets}\label{s6}
In this section we define an alpha limit set and (using such a set) introduce new types of minimality of CR-dynamical systems, all of them generalizing minimal dynamical systems. 

\begin{definition}
Let $(X,f)$ be a dynamical system and let  $\mathbf x\in \star_{i=1}^{\infty}\Gamma(f)^{-1}$. The set 
$$
\alpha_f(\mathbf x)=\{x\in X \ | \ \textup{ there is a subsequence of  the sequence } \mathbf x  \textup{ with limit } x\}
$$   
is called \emph{ \color{blue} the alpha limit set of $\mathbf x$}. 
\end{definition}

The following is a well-known result.  
\begin{theorem}
Let $(X,f)$ be a dynamical system. The following statements are equivalent.
\begin{enumerate}
\item \label{tritri3} $(X,f)$ is minimal.
\item \label{stiri4}
For each $x\in X$,  $T_{f}^{-}(x)\neq \emptyset$, and for each $\mathbf x\in \star_{i=1}^{\infty}\Gamma(f)^{-1}$,
$$
\alpha_f(\mathbf x)=X.
$$
\end{enumerate}
\end{theorem}
\begin{proof}
The proof is analogous to the proof of Theorem \ref{omomg}. We leave the details to a reader.
\end{proof}

\begin{definition}
Let $(X,G)$ be a CR-dynamical system, let $x_0\in X$ and let $\mathbf x\in T_G^{-}(x_0)$.  The set 
$$
\alpha_G(\mathbf x)=\{x\in X \ | \ \textup{ there is a subsequence of the sequence } \mathbf x \textup{ with limit } x\}
$$
is called \emph{ \color{blue}  the alpha limit set of $\mathbf x$} and we use  \emph{ \color{blue}  $\beta_G(x_0)$} to denote the set 
$$
\beta_G(x_0)=\bigcup_{\mathbf x\in T_G^{-}(x_0)}\alpha_G(\mathbf x).
$$
\end{definition}

\begin{definition}
Let $(X,G)$ be a CR-dynamical system.  We say that 
\begin{enumerate}
\item $(X,G)$ is \emph{ \color{blue} $1^{\alpha}$-minimal}, if for each $x\in X$,  $T_G^{-}(x)\neq \emptyset$, and for each $\mathbf x\in \star_{i=1}^{\infty}G^{-1}$,
$$
\alpha_G(\mathbf x)=X.
$$
\item $(X,G)$ is \emph{ \color{blue} $2^{\alpha}$-minimal}, if for each $x\in X$ there is $\mathbf x\in T_G^{-}(x)$ such that 
$$
\alpha_G(\mathbf x)=X.
$$
\item $(X,G)$ is \emph{ \color{blue} $3^{\alpha}$-minimal}, if for each $x\in X$,
$$
\beta_G(x)=X.
$$
\end{enumerate}
\end{definition}
\begin{observation}
Let $(X,G)$ be a CR-dynamical system. Then the following hold. 
\begin{enumerate}
\item If $(X,G)$ is $1^{\alpha}$-minimal, then $(X,G)$ is $2^{\alpha}$-minimal.
\item If $(X,G)$ is $2^{\alpha}$-minimal, then $(X,G)$ is $3^{\alpha}$-minimal.
\end{enumerate}
\end{observation}
Note that Example \ref{omegaPLUS} is also an example of a CR-dynamical system which is $2^{\alpha}$-minimal but is not $1^{\alpha}$-minimal. 
\begin{observation}\label{0666}
Let $(X,G)$ be a CR-dynamical system. Then the following hold. 
\begin{enumerate}
\item If $(X,G)$ is $1^{\alpha}$-minimal, if and only if $(X,G^{-1})$ is $1^{\omega}$-minimal.
\item If $(X,G)$ is $2^{\alpha}$-minimal, if and only if $(X,G^{-1})$ is $2^{\omega}$-minimal.
\item If $(X,G)$ is $3^{\alpha}$-minimal, if and only if $(X,G^{-1})$ is $3^{\omega}$-minimal.
\end{enumerate}
\end{observation}
\begin{theorem}\label{main2bb}
Let $(X,G)$ be a CR-dynamical system. Then the following hold.
\begin{enumerate}
\item \label{1trebb} $(X,G)$ is $1^{\alpha}$-minimal if and only if $(X,G)$ is $1^{\ominus}$-minimal. 
\item  \label{2trebb} $(X,G)$ is $2^{\alpha}$-minimal if and only if $(X,G)$ is $2^{\ominus}$-minimal.
\item  \label{3trebb} If $(X,G)$ is $3^{\alpha}$-minimal, then $(X,G)$ is $3^{\ominus}$-minimal.
\end{enumerate}
\end{theorem}
\begin{proof}
Let $(X,G)$ be a $1^{\alpha}$-minimal CR-dynamical system.  By Observation \ref{0666}, $(X,G)$ is  $1^{\alpha}$-minimal if and only if $(X,G^{-1})$ is $1^{\omega}$-minimal, and it follows from Theorem \ref{main666} that  $(X,G^{-1})$ is $1^{\omega}$-minimal if and only if $(X,G^{-1})$ is $1^{\oplus}$-minimal.  By Observation \ref{333},  $(X,G^{-1})$ is $1^{\oplus}$-minimal if and only if $(X,G)$ is $1^{\ominus}$-minimal. 

Let $(X,G)$ be a $2^{\alpha}$-minimal CR-dynamical system.  By Observation \ref{0666}, $(X,G)$ is  $2^{\alpha}$-minimal if and only if $(X,G^{-1})$ is $2^{\omega}$-minimal, and it follows from Theorem \ref{main666} that  $(X,G^{-1})$ is $2^{\omega}$-minimal if and only if $(X,G^{-1})$ is $2^{\oplus}$-minimal.  By Observation \ref{333},  $(X,G^{-1})$ is $2^{\oplus}$-minimal if and only if $(X,G)$ is $2^{\ominus}$-minimal. 

Let $(X,G)$ be a $3^{\alpha}$-minimal CR-dynamical system.  By Observation \ref{0666}, $(X,G)$ is  $3^{\alpha}$-minimal if and only if $(X,G^{-1})$ is $3^{\omega}$-minimal, and it follows from Theorem \ref{main666} that  $(X,G^{-1})$ is $3^{\oplus}$-minimal.  By Observation \ref{333},  $(X,G^{-1})$ is $3^{\oplus}$-minimal if and only if $(X,G)$ is $3^{\ominus}$-minimal. 
\end{proof}

\section{Preserving different types of minimality by topological conjugation}\label{s8}
The main results of this section are obtained in Theorem \ref{CMain}, where it is proved that any kind of minimality of a dynamical system is preserved by a topological conjugation. 
\begin{definition}
Let $X$ and $Y$ be metric spaces, and let $f:X\rightarrow X$ and $g:Y\rightarrow Y$ be  functions.  If there is a homeomorphism $\varphi:X\rightarrow Y$ such that 
$$
\varphi \circ f=g\circ \varphi,
$$
then we say that \emph{ \color{blue}  $f$ and $g$ are topological conjugates}. 
\end{definition}
The following is a well-known result.
\begin{theorem}
Let $(X,f)$ and $(Y,g)$ be dynamical systems. If $f$ and $g$ are topological conjugates, then
$$
(X,f) \textup{ is minimal } \Longleftrightarrow (Y,g)  \textup{ is minimal }.
$$
\end{theorem}
\begin{proof}
Let $\varphi:X\rightarrow Y$ be a homeomorphism such that 
$$
\varphi \circ f=g\circ \varphi.
$$
Suppose that $(X,f)$ is minimal and let $A$ be a non-empty closed subset of $Y$ such that $g(A)\subseteq A$.  Then $\varphi^{-1}(A)$ is a non-empty closed subset of $X$ such that 
$$
f(\varphi^{-1}(A))=\varphi^{-1}(g(A))\subseteq \varphi^{-1}(A).
$$
Therefore, $\varphi^{-1}(A)=X$ and $A=Y$ follows.  This proves that $(Y,g)$ is minimal. The proof of  the other implication is analogous.  
\end{proof}
The following definition generalizes the notion of topological conjugacy of continuous functions to the  topological conjugacy of closed relations. See \cite{BEK} for details.
\begin{definition}
Let $(X,G)$ and $(Y,H)$ be CR-dynamical systems. We say that \emph{ \color{blue}  $G$ and $H$ are topological conjugates} if there is a homeomorphism $\varphi:X\rightarrow Y$ such that for each $(x,y)\in X\times X$, the following holds
 $$
 (x,y)\in G  \Longleftrightarrow (\varphi(x), \varphi(y))\in H.
 $$
\end{definition}
In the rest of the paper,  we use $p_1,p_2:X\times X\rightarrow X$ to denote the projections $p_1(x,y)=x$ and $p_2(x,y)=y$ for all $(x,y)\in X\times X$, and $q_1,q_2:Y\times Y\rightarrow Y$ to denote the projections $q_1(x,y)=x$ and $q_2(x,y)=y$  for all $(x,y)\in Y\times Y$. Theorem \ref{CMain} is the main result of this section. We use the following lemmas in its proof.
\begin{lemma}\label{CLemma1}
Let $(X,G)$ and $(Y,H)$ be CR-dynamical systems and suppose that $G$ and $H$ are topological conjugates and let $\varphi:X\rightarrow Y$ be a homeomorphism such that for each $(x,y)\in X\times X$, 
 $$
 (x,y)\in G  \Longleftrightarrow (\varphi(x), \varphi(y))\in H.
 $$ 
Then the following hold.
\begin{enumerate}
\item $p_1(G)=X$ if and only if $q_1(H)=Y$.
\item $p_2(G)=X$ if and only if $q_2(H)=Y$.
\end{enumerate} 
\end{lemma}
\begin{proof}
Suppose that $p_1(G)=X$. To show that $q_1(H)=Y$, let $x\in Y$.  Let $z\in X$ such that $(\varphi^{-1}(x),z)\in G$ and let $y=\varphi(z)$. Then $(x,y)\in H$ and $x\in q_1(H)$ follows.  We have proved the implication from $p_1(G)=X$ to $q_1(H)=Y$. The proofs of the other three implications are analogous to the proof of this implication. We leave them to the reader. 
%
\end{proof}
\begin{lemma}\label{CLemma2}
Let $(X,G)$ and $(Y,H)$ be CR-dynamical systems and suppose that $G$ and $H$ are topological conjugates, let $\varphi:X\rightarrow Y$ be a homeomorphism such that for each $(x,y)\in X\times X$, 
 $$
 (x,y)\in G  \Longleftrightarrow (\varphi(x), \varphi(y))\in H,
 $$ 
 and let $A\subseteq X$. Then the following hold.
\begin{enumerate}
\item \label{last1} $A$ is $1$-invariant in $(X,G)$ if and only if $\varphi(A)$ is $1$-invariant in $(Y,H)$.
\item\label{last2} $A$ is $\infty$-invariant in $(X,G)$ if and only if $\varphi(A)$ is $\infty$-invariant in $(Y,H)$.
\item \label{last3} $A$ is $1$-backward invariant in $(X,G)$ if and only if $\varphi(A)$ is $1$-backward invariant in $(Y,H)$.
\item \label{last4} $A$ is $\infty$-backward invariant in $(X,G)$ if and only if $\varphi(A)$ is $\infty$-backward invariant in $(Y,H)$.
\end{enumerate} 
\end{lemma}
\begin{proof}
Suppose that $A$ is $1$-invariant in $(X,G)$. Obviously, since $\varphi:X\rightarrow Y$ is a homeomorphism and since $A$ is closed in $X$, also $\varphi(A)$ is closed in $Y$. Let $x\in \varphi(A)$ such that $x\in q_1(H)$ and let $z\in Y$ such that $(x,z)\in H$. Then $(\varphi^{-1}(x), \varphi^{-1}(z))\in G$ and, therefore, $\varphi^{-1}(x)\in A$ is such a point that $\varphi^{-1}(x)\in p_1(G)$.  Since $A$ is $1$-invariant in $(X,G)$, there is $w\in A$ such that $(\varphi^{-1}(x),w)\in G$. Fix such an element $w$ and let $y=\varphi(w)$.  Then $y\in \varphi(A)$ and $(x,y)\in H$. Therefore, $\varphi(A)$ is $1$-invariant in $(Y,H)$.  Next, suppose that $\varphi(A)$ is $1$-invariant in $(Y,H)$. Obviously, since $\varphi:X\rightarrow Y$ is a homeomorphism and since $\varphi(A)$ is closed in $Y$, also $A$ is closed in $X$. Let $x\in A$ such that $x\in p_1(G)$ and let $z\in X$ such that $(x,z)\in G$. Then $(\varphi(x), \varphi(z))\in H$ and, therefore, $\varphi(x)\in \varphi(A)$ is such a point that $\varphi(x)\in q_1(H)$.  Since $ \varphi(A)$ is $1$-invariant in $(Y,H)$, there is $w\in \varphi(A)$ such that $(\varphi(x),w)\in H$. Fix such an element $w$ and let $y=\varphi^{-1}(w)$.  Then $y\in A$ and $(x,y)\in G$. Therefore, $A$ is $1$-invariant in $(X,G)$.  This proves \ref{last1}.


The proofs of \ref{last2}, \ref{last3} and \ref{last4} are analogous to the proof of \ref{last1}.  We leave details to the reader. 
\end{proof}
\begin{lemma}\label{CLemma3}
Let $(X,G)$ and $(Y,H)$ be CR-dynamical systems and suppose that $G$ and $H$ are topological conjugates, let $\varphi:X\rightarrow Y$ be a homeomorphism such that for each $(x,y)\in X\times X$, 
 $$
 (x,y)\in G  \Longleftrightarrow (\varphi(x), \varphi(y))\in H.
 $$ 
Then the following hold.
\begin{enumerate}
\item \label{null} For each $x\in X$, $T^{+}_G(x)\neq \emptyset  $ if and only if $T^{+}_H(\varphi(x))\neq \emptyset$. 
\item  \label{eins} For each $\mathbf x\in \star_{i=1}^{\infty}G$,  $
\Cl(\Orbit(\mathbf x))=X $ if and only if $\Cl(\mathcal O_H^{\oplus}(\mathbf y))=Y
$
 where 
 $$
 \mathbf y=(\varphi(p_1(\mathbf x)),\varphi(p_2(\mathbf x)),\varphi(p_3(\mathbf x)),\ldots).
 $$
\item \label{zwei} For each $x\in X$, 
$
\Cl(\Orbitt(x))=X $ if and only if $\Cl(\mathcal U_H^{\oplus}(\varphi(x)))=Y
$.
\item  \label{drei} For each $\mathbf x\in \star_{i=1}^{\infty}G$, 
$
\omega_G(\mathbf x)=X $ if and only if $\omega_H(\mathbf y)=Y
$
 where 
 $$
 \mathbf y=(\varphi(p_1(\mathbf x)),\varphi(p_2(\mathbf x)),\varphi(p_3(\mathbf x)),\ldots).
 $$ 
\item \label{vier} For each $x\in X$,
$
\psi_G(x)=X $ if and only if $\psi_H(\varphi(x))=Y
$.
\item \label{null1} For each $x\in X$,
$
T^{-}_G(x)\neq \emptyset  $ if and only if $T^{-}_H(\varphi(x))\neq \emptyset
$.
\item  \label{eins1} For each $\mathbf x\in \star_{i=1}^{\infty}G^{-1}$, 
$
\Cl(\Orbitminus(\mathbf x))=X $ if and only if $\Cl(\mathcal O_H^{\ominus}(\mathbf y))=Y
$
 where 
 $$
 \mathbf y=(\varphi(p_1(\mathbf x)),\varphi(p_2(\mathbf x)),\varphi(p_3(\mathbf x)),\ldots).
 $$ 
\item \label{zwei1} For each $x\in X$,
$
\Cl(\mathcal U^{\ominus}_G(x))=X $ if and only if $\Cl(\mathcal U_H^{\ominus}(\varphi(x)))=Y
$.
\item  \label{drei1} For each $\mathbf x\in \star_{i=1}^{\infty}G$, 
$
\alpha_G(\mathbf x)=X $ if and only if $\alpha_H(\mathbf y)=Y
$
 where 
 $$
 \mathbf y=(\varphi(p_1(\mathbf x)),\varphi(p_2(\mathbf x)),\varphi(p_3(\mathbf x)),\ldots).
 $$ 
\item \label{vier1} For each $x\in X$, $\beta_G(x)=X $ if and only if $\beta_H(\varphi(x))=Y$.
\end{enumerate} 
\end{lemma}
\begin{proof}
First, we prove \ref{null}.  Let $x\in X$ such that $T^{+}_G(x)\neq \emptyset$ and let $\mathbf x\in T^{+}_G(x)$ and  $\mathbf y=(\varphi(p_1(\mathbf x)),\varphi(p_2(\mathbf x)),\varphi(p_3(\mathbf x)),\ldots)$.  Then  $\mathbf y\in T^{+}_H(\varphi(x))$ and, therefore, $T^{+}_H(\varphi(x))\neq \emptyset$.  To finish the proof of \ref{null}, let $x\in X$ such that $T^{+}_H(\varphi(x))\neq \emptyset$ and let $\mathbf y\in T^{+}_H(\varphi(x))$ and  $\mathbf x=(\varphi^{-1}(q_1(\mathbf y)),\varphi^{-1}(q_2(\mathbf y)),\varphi^{-1}(q_3(\mathbf y)),\ldots)$.  Then  $\mathbf x\in T^{+}_G(x)$ and, therefore, $T^{+}_G(x)\neq \emptyset$.  This completes the proof of \ref{null}.

Next, we prove \ref{eins}. Let $\mathbf x\in \star_{i=1}^{\infty}G$ and let $\mathbf y=(\varphi(p_1(\mathbf x)),\varphi(p_2(\mathbf x)),\varphi(p_3(\mathbf x)),\ldots)$.  First, suppose that $\Cl(\Orbit(\mathbf x))=X$.  To show that $\Cl(\Orbit(\mathbf y))=Y$, let $y\in Y$.  Then $\varphi^{-1}(y)\in \Cl(\Orbit(\mathbf x))$. Let $(x_n)$ be a sequence in $\Orbit(\mathbf x)$ such that $\displaystyle \lim_{n\to\infty}x_n=\varphi^{-1}(y)$.  Then $(\varphi(x_n))$ is a sequence in $\Orbit(\mathbf y)$ such that $\displaystyle \lim_{n\to\infty}\varphi(x_n)=y$. Therefore, $y\in \Cl(\Orbit(\mathbf y))$.   This proves the first implication of \ref{eins}.  To prove the other implication of \ref{eins}, suppose that $\Cl(\Orbit(\mathbf y))=Y$.  To show that $\Cl(\Orbit(\mathbf x))=X$, let $x\in X$.  Then $\varphi(x)\in \Cl(\Orbit(\mathbf y))$. Let $(y_n)$ be a sequence in $\Orbit(\mathbf y)$ such that $\displaystyle \lim_{n\to\infty}y_n=\varphi(x)$.  Then $(\varphi^{-1}(y_n))$ is a sequence in $\Orbit(\mathbf x)$ such that $\displaystyle \lim_{n\to\infty}\varphi^{-1}(y_n)=x$. Therefore, $x\in \Cl(\Orbit(\mathbf x))$.  This completes the proof of \ref{eins}.

The proofs of \ref{zwei},  \ref{drei},  \ref{vier},  \ref{null1},  \ref{eins1},  \ref{zwei1},  \ref{drei1},  and \ref{vier1} are straight forward and analogous to the proofs of \ref{null} and \ref{eins}. We leave them to the reader.
\end{proof}
\begin{theorem}\label{CMain}
Let $(X,G)$ and $(Y,H)$ be CR-dynamical systems and suppose that $G$ and $H$ are topological conjugates.  Then the following hold.
\begin{enumerate}
\item Let $k\in \{1,\infty,1^{\oplus},2^{\oplus},3^{\oplus},1^{\ominus},2^{\ominus},3^{\ominus},1^{\omega},2^{\omega},3^{\omega},1^{\alpha},2^{\alpha},3^{\alpha}\}$. 
 Then 
 $$
 (X,G) \textup{ is } k\textup{-minimal} \Longleftrightarrow  (Y,H) \textup{ is } k\textup{-minimal}.
 $$
 \item Let $k\in \{1,\infty\}$.   Then 
 $$
(X,G) \textup{ is } k\textup{-backward minimal} \Longleftrightarrow  (Y,H) \textup{ is } k\textup{-backward minimal}.
$$
 \end{enumerate}
\end{theorem}
\begin{proof}
Let $\varphi:X\rightarrow Y$ be a homeomorphism such that for each $(x,y)\in X\times X$, 
 $$
 (x,y)\in G  \Longleftrightarrow (\varphi(x), \varphi(y))\in H.
 $$
 We need to prove 16 statements. Their proofs are straight forward and  all of them  follow from Lemma \ref{CLemma2} or Lemma \ref{CLemma3}. We just give one of the proofs in details and leave the rest of the proofs to the reader.  
 
 We prove that
 $$
 (X,G) \textup{ is } 1\textup{-minimal} \Longleftrightarrow  (Y,H) \textup{ is } 1\textup{-minimal}.
 $$
 Suppose that $(X,G)$ is $1$-minimal. To prove that $(Y,H)$ is $1$-minimal, let $A$ be a non-empty closed subset of $Y$ which is $1$-invariant in $(Y,H)$.  By Lemma \ref{CLemma2},  $\varphi^{-1}(A)$ is a non-empty closed subset of $X$ which is $1$-invariant in $(X,G)$.  Since $(X,G)$ is $1$-minimal, it follows that $\varphi^{-1}(A)=X$. Therefore, $A=Y$.  It follows that $(Y,H)$ is $1$-minimal. 
 
 Suppose that $(Y,H)$ is $1$-minimal. To prove that $(X,G)$ is $1$-minimal, let $A$ be a non-empty closed subset of $X$ which is $1$-invariant in $(X,G)$.  By Lemma \ref{CLemma2},  $\varphi(A)$ is a non-empty closed subset of $Y$ which is $1$-invariant in $(Y,H)$.  Since $(Y,H)$ is $1$-minimal, it follows that $\varphi(A)=Y$. Therefore, $A=X$.  It follows that $(X,G)$ is $1$-minimal.

\end{proof}

%
%
%

\noindent I. Bani\v c\\
              (1) Faculty of Natural Sciences and Mathematics, University of Maribor, Koro\v{s}ka 160, SI-2000 Maribor,
   Slovenia; \\(2) Institute of Mathematics, Physics and Mechanics, Jadranska 19, SI-1000 Ljubljana, 
   Slovenia; \\(3) Andrej Maru\v si\v c Institute, University of Primorska, Muzejski trg 2, SI-6000 Koper,
   Slovenia\\
             {iztok.banic@um.si}           
     
				\-
				
		\noindent G.  Erceg\\
             Faculty of Science, University of Split, Rudera Bo\v skovi\' ca 33, Split,  Croatia\\
{{goran.erceg@pmfst.hr}       }    

                 \-
				
		\noindent R.  Gril Rogina\\
             Faculty of Natural Sciences and Mathematics, University of Maribor, Koro\v{s}ka 160, SI-2000 Maribor, Slovenia\\
{{rene.gril@student.um.si}       }    

                 	\-
					
  \noindent J.  Kennedy\\
             Lamar University, 200 Lucas Building, P.O. Box 10047, Beaumont, TX 77710 USA\\
{{kennedy9905@gmail.com}       }    

\end{document}